\documentclass[11pt]{article}
\usepackage{amssymb}
\usepackage{mathrsfs}
\usepackage{latexsym,amsmath,amssymb,amsfonts,epsfig,graphicx,cite,psfrag}
\usepackage{eepic,color,colordvi,amscd}
\usepackage{ebezier}
\usepackage{verbatim}

\usepackage{subfigure}

\newcommand{\pf}{\noindent {\bf Proof.} }

\newtheorem{theorem}{Theorem}[section]
\newtheorem{lemma}[theorem]{Lemma}
\newtheorem{coro}[theorem]{Corollary}

\newtheorem{ques}[theorem]{Question}

\def\qed{\hfill \rule{4pt}{7pt}}

\begin{document}

\title{Minimum co-degree condition for perfect matchings in $k$-partite $k$-graphs}

\author{Hongliang Lu\footnote{Partially supported by the National Natural
Science Foundation of China under grant No.11471257}\\
School of Mathematics and Statistics\\
Xi'an Jiaotong University\\
Xi'an, Shaanxi 710049, China\\
\medskip \\
Yan Wang\footnote{Partially supported by NSF grant DMS-1600738 through X. Yu} ~and  Xingxing Yu\footnote{Partially supported by NSF grant DMS-1600738
}\\
School of Mathematics\\
Georgia Institute of Technology\\
Atlanta, GA 30332}


\date{}

\maketitle

\begin{abstract}
Let $H$ be a  $k$-partite $k$-graph with $n$ vertices in each partition class,  and let
$\delta_{k-1}(H)$ denote the minimum co-degree of $H$.  
We characterize those $H$ with $\delta_{k-1}(H) \geq n/2$ 
and with no perfect matching. As a consequence we give an affirmative answer to the
following question of R\"odl and Ruci\'nski: If $k$ is even or $n \not\equiv 2 \pmod 4$, does $\delta_{k-1}(H) \geq n/2$ imply that $H$ has a perfect matching? We also give an example indicating that it
is not sufficient to impose this degree bound on only two types of $(k-1)$-sets.
\end{abstract}

\section{Introduction}

 A {\it hypergraph} $H$
consists of a vertex set $V(H)$ and an edge set $E(H)$ whose members
are subsets of $V(H)$. Let $H_1$ and $H_2$ be two hypergraphs. If   $V(H_1)\subseteq V(H_2)$ and $E(H_1)\subseteq E(H_2)$, then $H_1$ is called a \emph{subgraph} of $H_2$, denoted  $H_1\subseteq H_2$.
Let $k$ be a positive integer and $[k] := \{1,\ldots,k\}$.
For a set $S$, let ${S\choose k}:=\{T\subseteq S: |T|=k\}$.
A hypergraph $H$ is {\it $k$-uniform} if $E(H)\subseteq {V(H)\choose k}$, and a $k$-uniform hypergraph is also
called a {\it $k$-graph}. Given $T\subseteq V(H)$, let $H-T$ denote the subgraph of
$H$ with vertex set $V(H)-T$ and edge set $E(H-T)=\{e\in E(H)\ : \ e\subseteq V(H)-T\}$. 

Let $H$ be a $k$-graph and $S\in {V(H)\choose l}$ with $l \in [k]$.
The {\it neighborhood} of $S$ in $H$, denoted  $N_H(S)$, is the set of all $(k-l)$-subsets $U \subseteq V(H)$ such that $S \cup U \in E(H)$.
The {\it degree} of $S$ in $H$, denoted  $d_H(S)$,  is the size of $N_H(S)$.
For $l \in  [k]$, the {\it minimum $l$-degree} of $H$, denoted  $\delta_l(H)$, is the minimum degree over all $l$-subsets of $V(H)$.
Note that $\delta_{k-1}(H)$ is known as the minimum \emph{co-degree}  of $H$.

A \emph{matching} in a hypergraph $H$ is a subset of $E(H)$ consisting
of pairwise disjoint edges. A matching $M$ in a hypergraph $H$ is called a \emph{perfect matching}
if $V(M)=V(H)$.
R\" odl, Ruci\'nski, and Szemer\'edi \cite{RRS09} determined the minimum co-degree threshold function that ensures a perfect matching in a $k$-graph with $n$ vertices,
 for $n\equiv 0\pmod k$ and sufficiently large. This threshold function is $\frac{n}{2}-k+C$, where $C\in \{3/2,2,5/2,3\}$, depending on the parity of $n$ and $k$.
They \cite{RRS09} also proved that, for $n \not\equiv 0\pmod k$, the minimum co-degree threshold that ensures a matching $M$ in a $k$-graph $H$
with $|V(M)| \geq |V(H)| - k$ is between $\lfloor n/k\rfloor$ and $n/k+O(\log n)$, and  conjectured that this threshold function is $\lfloor n/k\rfloor$. This conjecture  was proved recently
by Han \cite{Han15}.
Treglown and Zhao \cite{TZ12, TZ13} determined the minimum $l$-degree threshold for perfect matchings in $k$-graphs for $k/2 \leq l \leq k-1$.

\medskip

A hypergraph $H$ is a {\it $k$-partite $k$-graph} with
\textit{partition classes} $V_1,\ldots,V_k$ if $V_1,\ldots,V_k$ is a
partition of $V(H)$ and $|e \cap V_i| = 1$ for all $e \in E(H)$ and $i \in [k]$.
We say that a set $S\subseteq V(H)$ is {\it legal} if $|S \cap V_i| \leq 1$ for $i \in [k]$.
For $l \in  [k]$,  the {\it minimum $l$-degree} of a $k$-partite
$k$-graph $H$, also denoted  $\delta_l(H)$, is the minimum degree over all legal $l$-subsets of $V(H)$.
Again, $\delta_{k-1}(H)$ is called the minimum \emph{co-degree}  of $H$.

 K\"uhn and Osthus \cite{Kuh06} showed that  the minimum co-degree
 threshold for the existence of a perfect matching in a $k$-partite $k$-graph with $n$
 vertices in each partition class is between $n/2$ and
 $n/2+\sqrt{2n\log n}$. Lu, Wang, and Yu \cite{Lu16} and, independently, Han, Zang, and Zhao \cite{HZZ16} showed that $n/k$ is the
minimum co-degree threshold for a $k$-partite $k$-graph $H$ with $n$ vertices in each partition class to admit a matching of size $|V(H)| - k$.

Aharoni, Georgakopoulos, and Spr\"ussel \cite{AGS09} obtained the
following stronger result:
Let $k \geq 3$ be a positive integer and
$H$ be a $k$-partite $k$-graph with
partition classes $V_1,\ldots,V_k$, each of size $n$. If $d_H(S) >
n/2$ for every legal
$(k-1)$-set $S$ contained in $V - V_1$, and if $d_H(T)\geq n/2$
for every legal $(k-1)$-set $T$ contained in $V-V_2$,
then $H$ has a perfect matching. Example 1 in  \cite{AGS09} (see the
graph $H_0(k,n)$ below) shows that this bound is best possible when $k$ is odd and $n\equiv 2 \pmod 4$.
Motivated by this result,
R\"odl and Ruci\'nski \cite{Rod09} asked the following
\begin{ques}
[R\"odl and Ruci\'nski \cite{Rod09}]
\label{rodl}
Let $k,n$ be integers with $k\geq 3$ and $n$ sufficiently large, and $H$ be a $k$-partite $k$-graph in which  each partition class has size $n$.
Assume that $k$ is even or $n \not\equiv 2 \pmod 4$. Is it true that if $\delta_{k-1}(H) \geq n/2$ then $H$ has a perfect matching?
If so, is it sufficient to impose this degree bound on only  two types of legal $(k-1)$-sets, similar to the above result of Aharoni, Georgakopoulos, and Spr\"ussel?
\end{ques}

Note that if $n$ is odd, it follows from the above result of Aharoni, Georgakopoulos, and Spr\"ussel that the answer to the
first part of Question~\ref{rodl} is affirmative.

We now describe an
example showing the tightness of the bound in Question~\ref{rodl}.
Let $k,n,d_i$, $i \in [k]$, be positive integers.
Let $H_0(d_1,\ldots,d_k;k,n)$ be a $k$-partite $k$-graph with
partition classes $V_1, \ldots , V_k$, and let   $D_i \subseteq V_i$
for $i \in [k]$, such that
  $|V_i| = n$ and  $|D_i| = d_i$ for $i \in [k]$, and
$E(H_0(d_1,\ldots,d_k;k,n))$ consists of those legal $k$-sets with
an even number of vertices (including zero) in $\bigcup_{i\in [k]}D_i$.
In particular,
we define
$H_0(k,n):=H_0(\lfloor n/2 \rfloor,\ldots,\lfloor n/2 \rfloor;k,n)$. 
When $k$ is odd and $n\equiv 2 \pmod 4$,  $H_0(k,n)$   is  Example 1 in \cite{AGS09}; in which case, $\delta_{k-1}(H_0(k,n))=
n/2$ and $H_0(k,n)$ admits no perfect matching  (as $\sum_{i\in [k]}|D_i|=kn/2$ is odd and every edge of $H_0(k,n)$ has an even number of vertices in $\bigcup_{i\in [k]}D_i$).

\medskip

{\bf Remark}. We point out that the answer to the second part of Question~\ref{rodl} is negative.
Let $k, n$ be positive integers such that $k$ is even or $n \equiv 0 \pmod 4$.
Let $J := H_0(n/2,n/2,\ldots,n/2,n/2+1;k,n)$ with partition classes $V_1, \ldots, V_k$ and let $D_i\subseteq V_i$ for $i\in [k]$ such that
$|D_i|=n/2$ for $i\in [k-1]$,  $|D_k| = n/2 + 1$,  and each edge
of $J$ has an even number of vertices in $\bigcup_{i\in [k]}D_i$.
Observe that all legal $(k-1)$-subsets of $V(J)$ intersecting $V_{k}$
have degree  at least $n/2$, and those legal $(k-1)$ sets contained in  $V(J)-V_k$ and intersecting $\cup_{i\in [k]}D_i$ an even number of times have degree $n/2-1$.
Moreover, $J$ has no perfect matching since $\sum_{i\in [k]}|D_i|= kn/2-1 \equiv 1 \pmod 2$ (as $k$ is even or $n \equiv 0 \pmod 4$).

\medskip

Our main result is the following, which implies an  affirmative answer to  the first part of
Question~\ref{rodl}.

\begin{theorem}
\label{main}
Let $k,n$ be integers with $k \geq 3$ and $n$ sufficiently large, and
let $H$ be a $k$-partite $k$-graph with $n$ vertices in each partition class.
Suppose $\delta_{k-1}(H) \geq \lfloor n/2\rfloor$.
Then $H$ has no perfect matching if, and only if,
\begin{itemize}
\item [ $(i)$] $k$ is odd, $n\equiv 2\pmod 4$,  and $H\cong H_0(k,n)$, or
\item [$(ii)$] $n$ is odd and there exist $d_i\in \{(n+1)/2,(n-1)/2\}$ for $i\in [k]$ such that   $\sum_{i=1}^k d_i$ is odd and $H \subseteq H_0(d_1,d_2,\ldots,d_k;k,n)$.
\end{itemize}
\end{theorem}

Our proof of Theorem~\ref{main} consists of two parts by considering whether or
not  $H$ is ``close'' to $H_0(k,n)$, which  is similar to arguments in \cite{Lo14, RRS09}.
Given two hypergraphs $H_1, H_2$ with $V(H_1)=V(H_2)$, let $c(H_1,H_2)$ be the minimum of $|E(H_1)\backslash E(H')|$
taken over all isomorphic copies $H'$ of $H_2$ with  $V(H') = V(H_2)$.
For a real number $\varepsilon > 0$,
we say that $H_2$ is \textit{$\varepsilon$-close} to $H_1$ if $V(H_1) = V(H_2)$ and $c(H_1,H_2)$ is less than $\varepsilon$ times the maximum possible number of edges on $V(H_2)$
(which is, for example, $\varepsilon n^k$ if $H_2$ is a $k$-partite $k$-graph with $n$ vertices in each partition class).

In Section 2, we deal with the case when  $H$ is $\varepsilon$-close to $H_0(k,n)$ for some sufficiently small $\varepsilon$.
In Section 3, we deal with the case  when  $H$ is not
$\varepsilon$-close to $H_0(k,n)$, using the absorbing method from \cite{RRS09} and
a recent result of the authors \cite{Lu16} (see Lemma~\ref{near_perfect}).

\section{Hypergraphs close to $H_0(k,n)$}

In this section, we prove Theorem~\ref{main}  for the case when $H$
is $\varepsilon$-close to $H_0(k,n)$ for some sufficiently small
$\varepsilon$.
Since we will be dealing with
$H_0(k,n)$, the following  notation for  ``even'' and ``odd'' degrees
(with respect to a given set $S$) will be convenient.
Let $H$ be a hypergraph.
 For $j \in \{0,1\}$, $v \in V(H)$, and $S\subseteq V(H)$, we define
$$d^{j}_{H,S}(v) := | \{ e \in E(H): v \in e \text{ and } |e \cap S| \equiv j \pmod 2 \}|.$$

\begin{lemma}
\label{near_critical}
Let $k \geq 3$ be a positive integer, and let   $\alpha, \varepsilon>0$ be small such
 that $\alpha<1/4$ and $\sqrt{\varepsilon}< \min\{1/(100 k^2),$ $1/(k (10k^2)^{k-1})\}$.
Then  for any $k$-partite $k$-graph $H$
with $n>100k^2$ vertices in each partition class, the following holds:
If $\delta_{k-1}(H) \geq (1/2-\alpha)n$, $H$ is $\varepsilon$-close to
$H_0(k,n)$, and $H\not\subseteq H_0(d_1,\ldots,d_k;k,n)$ for any $d_1,\ldots,d_k\in [\lceil (1/2-\alpha)n\rceil, \lfloor (1/2+\alpha)n\rfloor]$
with $\sum_{i=1}^k d_i$ odd, then  $H$ has a perfect matching.
\end{lemma}


\pf  Let  $H$ be a $k$-partite $k$-graph with $n$ vertices in each partition class such that
$\delta_{k-1}(H) \geq (1/2-\alpha)n$, $H$ is $\varepsilon$-close to $H_0(k,n)$, and
$H\not\subseteq H_0(d_1,\ldots,d_k;k,n)$ for any $d_1,\ldots,d_k\in [\lceil (1/2-\alpha)n\rceil, \lfloor (1/2+\alpha)n\rfloor]$
with $\sum_{i=1}^k d_i$ odd.
Let $$N:=\{v \in V(H): |N_{H_0(k,n)}(v) - N_H(v)| \geq
\sqrt{\varepsilon} n^{k-1}\}.$$
So each vertex in $N$ is contained in
at least $\sqrt{\varepsilon}n^{k-1}$ edges from $E(H_0(k,n))-E(H)$.
Note that  $$|N| \leq \sqrt{\varepsilon} k n;$$
for, otherwise,
\begin{eqnarray*}
|E(H_0(k,n))- E(H)|
&\geq& \frac{1}{k} \sum_{v \in N} |N_{H_0(k,n)}(v) - N_H(v)| \\
&>& \frac{1}{k}  |N|\sqrt{\varepsilon}n^{k-1}  \\
&>& \frac{1}{k} \sqrt{\varepsilon} kn \sqrt{\varepsilon}n^{k-1} \\
&=& \varepsilon n^k,
\end{eqnarray*}
contradicting the fact that $H$ is $\varepsilon$-close to $H_0(k,n)$.

\medskip

The rest of our proof is organized as follows. We first find a matching $M_1$ in $H$ that covers all  vertices
in $N$ (see Claim 2).
We then find a matching $M_2$  in $H-V(M_1)$
satisfying certain conditions (see
Claim 3).
Finally, we will show that there exists a
perfect matching in  $H-V(M_1)-V(M_2)$. The last part is easy when $k$
is even (see Claim 4), but needs more work when $k$ is odd (see Claims 5-8).

\medskip

To find a matching in $H$ that covers all vertices in $N$, we need to fix
some notation first.
For $i \in [k]$, let $B_i \subseteq V_i$ such that $|B_i| = \lfloor
n/2 \rfloor$
and each edge in $H_0(k,n)$ has an even number of vertices in $B := \cup_{j \in [k]} B_j$.
For $i \in [k]$, let $A_i := V_i - B_i$. The intuition for the notation below is that the vertices $v$ in $A_i\cap N$ (respectively, $B_i\cap N$)
with $ d^{0}_{H,B}(v) < n^{k-1}/8$ will be switched to $B_i'$ (respectively, $A_i'$).
For $i \in [k]$, let
$$A_i' := \left(A_i - \{v \in A_i \cap N: d^{0}_{H,B}(v) < n^{k-1}/8 \}\right) \cup \{v \in B_i \cap N: d^{0}_{H,B}(v) < n^{k-1}/8 \},$$
and
$$B_i' := \left(B_i - \{v \in B_i \cap N: d^{0}_{H,B}(v) < n^{k-1}/8 \}\right) \cup \{v \in A_i \cap N: d^{0}_{H,B}(v) < n^{k-1}/8 \}.$$
Let $A' := \cup_{j \in [k]} A_j'$ and $B' := \cup_{j \in [k]} B_j'$.

\medskip

Since $|N| \leq \sqrt{\varepsilon} k n$ and $|B_i|=\lfloor n/2\rfloor$, we
have $A_i'\neq \emptyset$ and $B_i'\neq \emptyset$ for $i\in [k]$ (as $n\ge 100k^2$).
In fact, for $i\in [k]$,
\begin{align}\label{Ai'_bound}
|A_i'| \geq |A_i| - |N| \geq (1/2 - \sqrt{\varepsilon}k)n
\end{align}
and
\begin{align}\label{Bi'_bound}
|B_i'| \geq |B_i| - |N| \geq (1/2 - \sqrt{\varepsilon}k)n - 1.
\end{align}
Moreover, for each $v\in V(H)$, the number of edges in $H$ containing
$v$ and intersecting $N-\{v\}$
is at most $|N|n^{k-2}$.

We now show that, for $v\in V(H)$,
\begin{align}\label{degree_bound_H-N}
d_{H-(N-\{v\}),B'}^{0}(v) \ge (1/8 - \sqrt{\varepsilon}k ) n^{k-1}.
\end{align}
 First, suppose
 $v \in (A \cap A') \cup (B \cap B')$. Then
 $B'-(N-\{v\})=B-(N-\{v\})$, and $d_{H,B}^0(v) \geq  n^{k-1}/8$ by
 definition of $A',B'$.  So $$d_{H-(N-\{v\}),B-(N-\{v\})}^0(v) \geq d_{H,B}^{0}(v) - |N|n^{k-2} \geq  n^{k-1}/8 - \sqrt{\varepsilon} k n^{k-1}.$$
Hence,
\begin{align*}
d_{H-(N-\{v\}),B'}^{0}(v)
&= d_{H-(N-\{v\}),B'-(N-\{v\})}^{0}(v) \\
&= d_{H-(N-\{v\}),B-(N-\{v\})}^{0}(v) \\
&\ge (1/8 - \sqrt{\varepsilon}k ) n^{k-1}.
\end{align*}
Now assume $v \in (A \cap B')\cup (A'\cap B)$. Then  $v\in N$ and
$d_{H,B}^0(v) < n^{k-1}/8$.  So $B'-(N-\{v\})=B-(N-\{v\})-\{v\}$,
Since $\delta_{k-1}(H) \geq (1/2 - \alpha)n$, it follows that $d_{H}(v) \geq (1/2 - \alpha)n^{k-1}$. Thus, since $n\ge 100k^2$ and $\alpha<1/4$,
$$d_{H,B}^1(v) \geq (1/2 - \alpha)n^{k-1} - d_{H,B}^0(v) > (1/2 - \alpha)n^{k-1} - n^{k-1}/8 > n^{k-1}/8.$$
Therefore,  $d_{H-(N-\{v\}),B-(N-\{v\}) }^{1}(v) \geq d_{H,B}^{1}(v) - |N|n^{k-2} \geq (1/8 - \sqrt{\varepsilon}k ) n^{k-1}$.
Hence,
\begin{align*}
d_{H-(N-\{v\}),B'}^{0}(v)
&= d_{H-(N-\{v\}),B'-(N-\{v\})}^{0}(v) \\
&= d_{H-(N-\{v\}),B-(N-\{v\}) }^{1}(v) \\
&\ge (1/8 - \sqrt{\varepsilon}k ) n^{k-1}.
\end{align*}

We now begin our process of finding matchings $M_1$ and $M_2$. First, we need to make $|B'|$ even.
\medskip

\textit{Claim 1.}  Either  $|B'|$ is even (in which case let
$e_0=\emptyset$; so $|B'-e_0|$
is even), or there exists an edge $e_0\in E(H)$ such that
 $|B' - e_0|$ is even.

We may assume that $|B'|$ is odd and $|B'-e|$ is odd for every  $e \in
E(H)$;  as, otherwise, Claim 1 holds.
Then $|B'\cap e|$ is even for all $e\in E(H)$.
Hence $H \subseteq H_0(d_1,\ldots,d_k;k,n)$, where $d_i = |B_i'|$ for $i \in [k]$, $\sum_{i\in [k]}d_i=|B'|$ is odd,
and $B_1',\ldots, B_k'$ play the roles of $D_1,\ldots, D_k$, respectively, in the definition of $H_0(k,n)$.

Let $v_i \in A_i'$ and $u_i \in B_i'$ for $i \in [k]$, and
let $S := \{v_1, \ldots, v_k\}$.
Then for $i \in [k]$, since  $|B'\cap e|$ is even for all $e\in E(H)$,
we have
$$n - d_i = |A_i'| \geq d_H(S - \{v_i\}) \geq \delta_{k-1}(H) \geq (1/2 - \alpha)n; $$
so $d_i \leq \lfloor (1/2 + \alpha)n\rfloor.$ Moreover, for $i \in
[k]$, let $j \in [k] - \{i\}$. Again, since  $|B'\cap e|$ is even for
all $e\in E(H)$, we have
$$d_i = |B_i'| \geq d_H( (S \cup \{u_j\}) - \{v_i, v_j\}) \geq \delta_{k-1}(H) \geq (1/2 - \alpha)n; $$
so $d_i \geq \lceil(1/2- \alpha)n\rceil$. This contradicts the assumption that $H \not\subseteq H_0(d_1,\ldots,d_k;k,n)$ for any
$d_1,\ldots,d_k\in [\lceil (1/2-\alpha)n\rceil, \lfloor (1/2+\alpha)n\rfloor]$ with $\sum_{i=1}^k d_i$ odd.  $\Box $

\medskip

Note that for each  $v\in N - e_0$, the number of edges in $H$
containing $v$ and a vertex of $e_0$ is at most $kn^{k-2}$.
Thus by (\ref{degree_bound_H-N}), we have
\begin{align}\label{degree_bound_H-N-e0}
d_{(H-e_0)-(N-\{v\}),B'-e_0}^{0}(v)\ge (1/8- 2\sqrt{\varepsilon}k ) n^{k-1} - k n^{k-2} > n^{k-1}/10,
\end{align}
where the last inequality holds since $\sqrt{\varepsilon} < 1/(100k^2)$ and $n\ge 100k^2$.

\medskip

\textit{Claim 2. } There exists a matching $M_1$ in $H-e_0$ such that
\begin{itemize}
\item [$(i)$] $|M_1|\le \sqrt{\varepsilon} kn$,
\item [$(ii)$] $N-e_0\subseteq V(M_1)$, and
\item [$(iii)$]  $|e \cap (B'-e_0)| \equiv 0 \pmod 2$ for all $e\in M_1$.
\end{itemize}

Let $M_1:= \emptyset$ if $N - e_0 = \emptyset$.
Now assume  $N - e_0 \neq \emptyset$, 
and we construct $M_1$ by matching vertices in $N$ greedily.
Let $v_1 \in N - \{e_0\}$.
Since $d_{(H-e_0)-(N-\{v_1\}),B'-e_0}^{0}(v_1) >n^{k-1}/10$ (by (\ref{degree_bound_H-N-e0}))
and $n\ge 100k^2$, there exists an edge  $e_1$ in $H-e_0$, such
that $v_1 \in e_1$ and $|e_1 \cap (B'-e_0)| \equiv 0 \pmod 2$.

Now  suppose
we have found a
matching $\{e_1,\ldots,e_t\}$ in $H-e_0$ for some $t\ge 1$,
such that, for $i\in [t]$,
$e_i \cap (N - e_0) \neq \emptyset$ and
$|e_i \cap (B' - e_0)| \equiv 0 \pmod 2$.
If $N - e_0 \subseteq \cup_{i \in [t]} e_i$, then $M_1 := \{e_1,
\ldots, e_t\}$ is the desired matching (as $t<|N|\le \sqrt{\varepsilon}kn$).
So let $v_{t+1} \in N - e_0$ and $v_{t+1} \not \in \cup_{i \in [t]} e_i$.
Note that $t<|N|\le \sqrt{\varepsilon}kn$ and that the number of
edges in $H-e_0$ containing $v_{t+1}$ and  a vertex
from $\cup_{i\in   [t]} e_i$ is
at most $$tkn^{k-2} \leq \sqrt{\varepsilon} k^2 n^{k-1} \leq n^{k-1}/100$$
as  $\sqrt{\varepsilon} < 1/(100k^2$.
Since $d_{(H-e_0)-(N-\{v_{t+1}\}),B'-e_0}^{0}(v_{t+1}) >n^{k-1}/10$ (by (\ref{degree_bound_H-N-e0})),
there exists $e_{t+1}$ in $(H-e_0) - \cup_{i\in   [t]} e_i$ such that $v_{t+1}\in e_{t+1}$
and $|e_{t+1} \cap (B' - e_0)| \equiv 0 \pmod 2$.

Therefore, continuing this process (at most $|N-e_0|$ steps), we obtain the desired matching for
Claim 2. $\Box $

\medskip

Let $H' := (H -e_0)- V(M_1)$.
For $i\in[k]$, let  $C_i := A_i' - (V(M_1) \cup e_0)$, $D_i := B_i' - (V(M_1) \cup e_0)$ and $D := \cup_{i \in [k]} D_i$.
By Claim 2, $N\cap V(H')=\emptyset$; so  for $i \in [k]$,
\begin{align}\label{D_i_in_B_i}
D_i \subseteq B_i.
\end{align}
Note that $|D|$ is even (by Claims 1 and 2).
Since $|M_1| \leq \sqrt{\varepsilon} k n$, it follows from (\ref{Ai'_bound}) and (\ref{Bi'_bound}) that for $i\in [k]$,
\begin{align}\label{C_i_bound}
|C_i| \geq |A_i'| - (|M_1| + 1) \geq (1/2 - \sqrt{\varepsilon}k)n  -  (\sqrt{\varepsilon} k n + 1)
\geq (1/2 - 2\sqrt{\varepsilon}k)n - 1
\end{align}
and
\begin{align}\label{D_i_bound}
|D_i| \geq |B_i'| - (|M_1| + 1) \geq ((1/2 - \sqrt{\varepsilon}k)n - 1) -  (\sqrt{\varepsilon} k n + 1)
\geq (1/2 - 2\sqrt{\varepsilon}k)n - 2.
\end{align}

\medskip

\textit{Claim 3. } There exists a matching $M_2$ in $H'$ such that
\begin{itemize}
\item [$(i)$] $|M_2| \leq 8 \sqrt{\varepsilon}k^2 n$,
\item [$(ii)$] $|D_i - V(M_2)|=|D_1-V(M_2)|$  for $i \in [k] $, and
\item [$(iii)$] $|D - V(M_2)|$ is even.
\end{itemize}

Without loss of generality, we may assume that $|D_1| \geq |D_2| \geq
\cdots \geq |D_k|$. If $|D_1|=|D_k|$ then $M_2=\emptyset$ gives the
desired matching for Claim 3. So assume $|D_1| - |D_k| > 0$.
We construct an auxiliary graph and use a perfect matching in this graph to find $M_2$.

Let $r \in \{0,1\}$ such that $|D_1|+r$ is even.
Let $G$ be the complete $k$-partite 2-graph and let  $W_1,...,W_k$ be
the partition classes of $G$, such that $|W_i| = (|D_i|-|D_k|)+ (|D_1| +r)- |D_k|$ for $i \in [k]$.
Then
$|W_1| \geq |W_2| \geq \ldots \geq |W_k|$ and
$$|V(G)| = \sum_{i \in [k]} |W_i|
= \left(\sum_{i \in [k]} |D_i|\right) + k  ( |D_1| + r ) - 2k |D_k|.$$
Since $\sum_{i \in [k]} |D_i|$ and $|D_1| + r$ are even, $|V(G)|$ is also even.

We now use Tutte's 1-factor theorem to show that $G$ has a
perfect matching. For $S \subseteq V(G)$, let $o(G-S)$ denote the
number of connected components of  $G - S$ of odd order.
If $S = \emptyset$, then $o(G-S) = 0 \leq |S|$. Now assume $S\ne \emptyset$.
Since $G$ is a complete $k$-partite 2-graph and $|W_1| \geq |W_2| \geq \ldots \geq |W_k|$,  if $1 \leq |S| < \sum_{i
  \in [k] - \{1\}} |W_i|$ then $o(G-S) \leq 1 \leq |S|$, and if $|S| \geq \sum_{i \in [k] - \{1\}} |W_i| \geq (k - 1) (|D_1|-|D_k| + r)$
then $o(G-S) \leq |W_1| \leq 2 |D_1|-2|D_k| + r \leq |S|$ (as $k \geq 3$).
Thus, by Tutte's 1-factor theorem, $G$ has a perfect matching,
say  $T$.

Since $|C_i| \geq (1/2 - 2\sqrt{\varepsilon}k)n - 1$ (by (\ref{C_i_bound})),
$|D_i| \leq n  - |C_i| \leq (1/2 + 2\sqrt{\varepsilon}k)n + 1$.
So by (\ref{D_i_bound}), $|D_1|-|D_k| \leq 4 \sqrt{\varepsilon}k n + 3.$
Hence,
$$|T| = |V(G)|/2 = \left(\sum_{i \in [k]}|W_i| \right)/2 \leq (k/2) (2|D_1|-2|D_k| + 1) \leq 8 \sqrt{\varepsilon}k^2 n.$$

Let $T=\{f_1,f_2,...,f_{|T|}\}$. Corresponding to each $f_i$ we find an
edge $g_i$ of $H'$ such that $\{g_1, \ldots, g_{|T|}\}$ gives
the desired matching $M_2$ for Claim 3.

Let $g_0=\emptyset$ and we find $g_1, \ldots, g_{|T|}$ in order.
Suppose we have found $g_t$ for some $t$,
with $0\le t\le |T|-1$. We describe how to find $g_{t+1}$ using $f_{t+1}$.
Let $f_{t+1} \subseteq W_p \cup W_q$, where $p,q \in [k]$.
By (\ref{C_i_bound}) and (\ref{D_i_bound}),
$\min\{|C_j|,|D_j|\}\ge (1/2 - 2\sqrt{\varepsilon}k)n - 2$ for $j\in [k]$. Then,
since $|T|\le 8 \sqrt{\varepsilon}k^2 n$, $n\ge 100k^2$, and $\sqrt{\varepsilon}<1/(k(10k^2)^{k-1})$, we have, for $j\in [k]$,
$$|C_j-\cup_{i \in [t]} g_i|>n/10 \mbox{ and } |D_j-\cup_{i \in [t]} g_i|>n/10.$$
So let $v_p \in D_p - \cup_{i \in [t]} g_i$. There
exist $v_q \in D_q - \cup_{i \in [t]} g_i$ and $v_j \in C_j - \cup_{i \in [t]} g_i$ for $j \in [k] - \{p,q\}$
such that
$g_{t+1} := \{v_1,...,v_k\} \in E(H')$; for,  otherwise,
$$ |N_{H_0(k,n)}(v_p) - N_H(v_p)|> ( n/10 )^{k-1} > \sqrt{\varepsilon} n^{k-1},$$
as $\sqrt{\varepsilon}<1/(k(10k^2)^{k-1}$, contradicting the fact that $v_p \not\in N$.
Clearly, $|g_{t+1}\cap D|=2$.

Therefore, $M_2 := \{g_1,...,g_{|T|}\}$  is a matching in $H'$ such
that, for $i \in [k]$, $$|D_i - V(M_2)| =|D_i|-|W_i|= 2|D_k| - |D_1| - r > 0,$$
where the inequality holds because of (\ref{D_i_bound}),  $|D_1|-|D_k|\le 4\sqrt{\varepsilon} kn+3$, $\sqrt{\varepsilon}<1/(100k^2)$, and $n\ge 100k^2$.
Moreover,  $|g_j \cap D| = 2$ for $j\in [|T|]$.
Hence,  since $|D|$ is even (by Claims 1 and 2), $|D - V(M_2)|$ is even.
$\Box$

\medskip

Let $H'' := H' - V(M_2)$ and, for $i\in [k]$, let  $D_i' :=  D_i -
V(M_2)$ and $C_i' := C_i - V(M_2)$.
Let $D' := \cup_{i \in [k]} D_i'$ and $C' := \cup_{i \in [k]}
C_i'$. Note that $|D'|$ is even, as $|D-V(M_2)|$ is even (by Claim 3).
Since $|M_2| \leq 8 \sqrt{\varepsilon} k^2 n$ (by Claim 3),
it follows from (\ref{C_i_bound}) and (\ref{D_i_bound}) that, for $i\in [k]$,
\begin{align}\label{C_i'_D_i'_bound}
\min\{|C_i'|,|D_i'|\} =\min \{|C_i|,|D_i|\} - |M_2| \geq (1/2 - 2\sqrt{\varepsilon}k)n - 2 - 8 \sqrt{\varepsilon} k^2 n.
\end{align}

\medskip

\textit{Claim 4.} We may assume that $k$ is odd.

For,  suppose $k$ is even. We show that both $H'' - C'$ and $H'' - D'$
have perfect matchings;  hence the assertion of the lemma
holds. Below, we only show that $H'' - C'$ has a perfect matching,
since the argument for $H''-D'$ is the same (by substituting  (\ref{C_i_bound}) for (\ref{D_i_bound}) and by
exchanging the roles of $C_i'$ and $D_i'$).

Let $M$ be a maximum matching in $H'' - C'$. Then $(H''-C')-V(M)=H[D' - V(M)]$ has no edge.
We claim that $|M| \geq n/4$. For, otherwise, $D_1'-V(M)\ne
\emptyset$ by (\ref{D_i_bound})
 (as $\sqrt{\varepsilon}<1/(100k^2)$ and $n\ge 100k^2$).
Let  $v \in D_1' - V(M)$. Since
$k$ is even and $D_i' \subseteq B_i$ for $i \in [k]$ (by (\ref{D_i_in_B_i})), and because $H[D'-V(M)]$ has no edge,
we have
\begin{eqnarray*}
|N_{H_0(k,n)}(v) - N_H(v)|
&\geq &|D_2' - V(M)|  |D_3' - V(M)| \cdots |D_k' - V(M)| \\
&\geq & ( (1/2 - 2\sqrt{\varepsilon}k)n - 2 - 8 \sqrt{\varepsilon} k^2 n - n/4 )^{k-1} \quad \mbox{  (by (\ref{C_i'_D_i'_bound}))} \\
&> & ( n/10 )^{k-1} \quad \mbox{ (since $\sqrt{\varepsilon} < 1/(100k^2)$ and $n\ge 100k^2$)}\\
&> & \sqrt{\varepsilon} n^{k-1} \quad \mbox{ (since $\sqrt{\varepsilon} < 1/(k(10k^2)^{k-1})$},
\end{eqnarray*}
contradicting the fact that $v \not\in N$.

Now, suppose for a contradiction, that $M$ is not a perfect matching in $H'' - C'$.
Then there exists $u_i \in D_i' - V(M)$ for $i \in [k]$.
Note that $|M| \geq n/4 > k-1$ (as $n\ge 100k^2$).

Let $\{e_1,\ldots,e_{k-1}\}$ be an arbitrary $(k-1)$-subset of $M$, and write
$e_i := \{v_{i,1},\ldots,v_{i,k}\}$ with $v_{i,j} \in D_j'$ for $i \in [k-1]$ and  $j \in [k]$.
For $j \in [k]$, let $f_j := \{u_j, v_{1,j+1},
v_{2,j+2},\ldots,$ $v_{k-1,j+k-1}\}$, with the addition in the subscripts
modulo $k$ (except we write $k$ for $0$). Note that $f_1, \ldots, f_k$ are pairwise disjoint.
Since $D_i' \subseteq B_i$ for $i \in [k]$ (by (\ref{D_i_in_B_i})),
and $k$ is assumed to be even, it follows that  $f_j\in E(H_0(k,n))$ for
$j\in [k]$.

If $f_i \in E(H'')$ for all $i \in [k]$ then $M':= (M \cup
\{f_1,\ldots,f_k\}) - \{e_1,\ldots,e_{k-1}\}$ is a matching in $H$
and  $|M'| = |M| + 1 > |M|$, contradicting the maximality of $|M|$.

Hence, $f_j\not\in E(H)$ for some $j \in [k]$.
Note that there are $\binom{|M|}{k-1}$ choices of $\{e_1,\ldots, e_{k-1}\} \subseteq M$.
Hence,
\begin{eqnarray*}
& & |\{e \in E(H_0(k,n)) - E(H): |e\cap \{u_i: i\in [k]\}|=1\}|\\
&\geq & |M|(|M| - 1) \cdots (|M| - (k - 1) + 1) / (k-1)! \\
&> & {(n/4 - k + 2)}^{k-1} / (k-1)! \\
&> & (n/10k)^{k-1} \quad \mbox{ (since $n\ge 100k^2$)}\\
&> & k  \sqrt{\varepsilon} n^{k-1}  \quad \mbox{ (since $\sqrt{\varepsilon} < 1/(k(10k^2)^{k-1})$}.
\end{eqnarray*}
This implies that there exists $i \in [k]$ such that
$|N_{H_0(k,n)}(u_i) - N_{H}(u_i)| > \sqrt{\varepsilon} n^{k-1}$,
contradicting the fact that $u_i \not\in N$. $\Box$

\medskip

Next claim guarantees a divisibility condition for $|D'|$, which will be used in the proof of Claim 7.

\medskip

\textit{Claim 5.} There exists a matching $M_3$ in $H''$ such that
\begin{itemize}
\item [$(i)$]  $|M_3|\leq k^2/2$, and
\item [$(ii)$] $|D_i' - V(M_3)|=|D_1'-V(M_3)| \equiv 0 \pmod {k-1}$ for $i \in [k]$.
\end{itemize}

Let $0 \leq s \leq k-2$ be such that $|D_1'| \equiv s \pmod {k-1}$.
We may assume that $s \neq 0$; for, otherwise, $M_3 = \emptyset$ gives
the desired matching for Claim 5.
Moreover, since $k$ is odd (by Claim 4) and  $|D'|=k|D_1'|$ is even,
it follows that $s$ is even.

We now construct $M_3$, starting with the empty matching
$T_0=\emptyset$. Suppose for some $j \in [s/2]$,
we have constructed a matching $T_{j-1}$ in $H''$ with $|T_{j-1}| = k(j-1)$.
Since $\sqrt{\varepsilon}<1/(100k^2)$ and $n\ge 100k^2$, it follows from (\ref{C_i'_D_i'_bound}) that, for  $i \in [k]$,
$$\min\{|C_i'-V(T_{j-1})|, |D_i'-V(T_{j-1})|\}\ge (1/2-2\sqrt{\varepsilon}k)n-2 - 8 \sqrt{\varepsilon} k^2 n-k(j-1) > 0.$$

For $i\in [k]$,  let $v_{j,i}\in D_i' - V(T_{j-1})$.
We claim that there exist $v_{j,i+1} \in D_{i+1}' - V(T_{j-1})$ and
$u_{j,l} \in C_{l}' - V(T_{j-1})$ for $l \in [k] - \{i, i+1\}$, such
that $e_{j,i} := \{ v_{j,i}, v_{j,i+1}, u_{j,l} : l \in [k] - \{i,
i+1\}  \} \in E(H'')$ (with addition in the subscripts  modulo $k$ except we use $k$ for $0$) and $\{e_{j,i}:i\in [k]\}$ is a matching in $H''$.
For, otherwise, since $D_i' \subseteq B_i$ by ((\ref{D_i_in_B_i})), we have
\begin{eqnarray*}
|N_{H_0(k,n)}(v_{j,i}) - N_H(v_{j,i})|
&\geq & |D_{i+1}' - V(T_{j-1})-\cup_{i\in [k]}e_{j,i}|  \prod_{l \in [k] - \{i, i+1\}} |C_l' - V(T_{j-1})| \\
&\geq & ( (1/2 - 2\sqrt{\varepsilon}k)n - 2 - 8 \sqrt{\varepsilon} k^2 n - k^2/2 )^{k-1} \\
&> & ( n/10 )^{k-1}  \mbox{ (since $\sqrt{\varepsilon} < 1/(100k^2)$ and $n\ge 100k^2$)}\\
&> & \sqrt{\varepsilon} n^{k-1}  \mbox{ (since $\sqrt{\varepsilon} < 1/(k(10k^2)^{k-1})$},
\end{eqnarray*}
contradicting the fact that $v_{j,i} \not\in N$.

Let $T_j = T_{j-1} \cup \{e_{j,i} : i \in [k]\}$. Then $T_j$ is a
matching in $H''$ for $j\in [s/2]$.
Let $M_3 := T_{s/2} = \{e_{j,i}: j \in [s/2] \mbox{ and } xs i \in [k] \}$. Then $|M_3|\le
k^2/2$. Note that, for $i\in [k]$, the edges in $T_j-T_{j-1}$ uses
exactly two vertices of $D_i'$. Thus, for $i\in [k]$,  $|D_i' - V(M_3)|=|D_1' - V(M_3)| = |D_1'| - s \equiv 0
\pmod {k-1}$. $\Box$

\medskip

Let $H^* := H'' - V(M_3)$ and, for $i\in [k]$, let $D_i^* := D_i' - V(M_3)$ and $C^*_i := C_i' - V(M_3)$.
Let $D^* := \cup_{i \in [k]} D_i^*$ and $C^*:= \cup_{i \in [k]} C_i^*$.
Since $|M_3| \leq k^2/2$ (by Claim 5), it follows from (\ref{C_i'_D_i'_bound}) that
\begin{align}\label{C_i*_D_i*_bound}
\min\{|C_i^*|,|D_i^*|\}\ge \min\{|C_i'|,|D_i'|\}-|M_3|\ge (1/2 -
2\sqrt{\varepsilon}k)n - 2 - 8 \sqrt{\varepsilon} k^2 n - k^2/2.
\end{align}
By Claim 5, $|D_i^*| =|D_1^*|\equiv 0 \pmod {k-1}$ for $i \in [k]$.

We will show that  $H^*$ has a perfect matching  using
edges of special types. For any $e\in
E(H^*)$, if $e \subseteq C^*$ then we say that $e$ is of
\textit{0-type}, and if
$|e \cap C^*| =|e\cap C_j^*|= 1$  for some $j \in [k]$ then we say that $e$ is of
\textit{j-type}.
For convenience, let $$\tau := 1/(9 k).$$

\medskip

\textit{Claim 6.}
$H^*$ has pairwise disjoint matchings  $M^0, M^1,
\ldots, M^k$, such that for $i\in [k]\cup \{0\}$,
\begin{itemize}
\item [$(i)$] $|M^i|=\lfloor \tau n\rfloor$, and
\item [$(ii)$] each edge in $M^i$ is  of $i$-type.
\end{itemize}

We construct $M^0,M^1,\ldots,M^k$ in the order listed.
Let $T^0$ be a matching in $H^*$ such that $V(T^0)\subseteq C^*$ and,
subject to this, $|T^0|$ is maximum. Then $C^* - V(T^0)$ has no edge.  We claim that $|T^0| \geq \lfloor \tau n\rfloor$; for, otherwise,
 $|C_i^* - V(T^0)|\ge |C_i^*|-\tau n$ for $i\in [k]$ and, hence, for any $v \in C_1^* - V(T^0)$,
\begin{eqnarray*}
|N_{H_0(k,n)}(v) - N_H(v)|
&\geq & |C_2^* - V(T^0)| |C_3^* - V(T^0)| \cdots |C_k^* - V(T^0)| \\
&\geq & ( (1/2 - 2\sqrt{\varepsilon}k)n - 2 - 8 \sqrt{\varepsilon} k^2 n - k^2/2 - \tau n )^{k-1}  \mbox{ (by (\ref{C_i*_D_i*_bound}))} \\
&> & ( n/10 )^{k-1} \mbox{ (since $\sqrt{\varepsilon} < 1/(100k^2)$ and $n\ge 100k^2$)}\\
&> & \sqrt{\varepsilon} n^{k-1} \mbox{ (since $\sqrt{\varepsilon} < 1/(k(10k^2)^{k-1})$},
\end{eqnarray*}
contradicting the fact that $v \not\in N$.
Let $M^0$ be a set of any $\lfloor \tau n\rfloor$ edges in $T^0$.

Now suppose for some $j\in [k]$,   we have found matchings
$M^0,M^1,\ldots,$ $M^{j-1}$ in $H^*$ such that $M^i$ (for $i=0,
\ldots, j-1$) consists of $\lfloor \tau n\rfloor$ edges of $i$-type.
Let $T^j$ be a matching in $H^* - \cup_{i=0}^{j-1} V(M_i)$ such that each edge in $T_j$ is of $j$-type and, subject to this, $|T^j|$ is maximum.

We claim that $|T^j| \geq \lfloor \tau n\rfloor$. For, suppose $|T^j|<
\lfloor \tau n\rfloor$.
Then, since $C_j^* \cap V(M^i) = \emptyset$ for $i \in [j-1]$ and $|V(M^0)|= \lfloor \tau n\rfloor$, it follows from  (\ref{C_i*_D_i*_bound}) that
$$|C_j^* - V(M^0 \cup M^1 \cup \ldots \cup  M^{j-1}) - V(T^j)| > (1/2 - 2\sqrt{\varepsilon}k)n - 2 - 8 \sqrt{\varepsilon} k^2 n - k^2/2 - k\tau n > 0,$$
where the second inequality holds because $\tau =1/(9k)$,
$\sqrt{\varepsilon}<1/(100k^2)$, and $n\ge 100k^2$.
.
So let $v$ be a vertex in $C_j^* - V(M^0 \cup M^1 \cup \ldots \cup  M^{j-1}) - V(T^j)$.
We claim that there exists an edge $f$ of $j$-type in
$H^*-V(M^0\cup M^1 \cup \ldots \cup M^{j-1}) - V(T^j)$ with  $v \in f$;
as, otherwise, since $D_i^* \subseteq B_i$ for $i \in [k]$ (by (\ref{D_i_in_B_i})) and $k$ is odd,
\begin{eqnarray*}
| N_{H_0(k,n)}(v) - N_H(v) |
&\geq & \prod_{l \in [k]-\{j\}} \big| D_l^* -  V(M^0\cup M^1 \cup \ldots \cup M^{j-1}) - V(T^j) \big| \\
& \geq & ( ( 1/2 - 2\sqrt{\varepsilon}k)n - 2 - 8 \sqrt{\varepsilon}
         k^2 n - k^2/2 - k \tau n )^{k-1} \mbox{ (by (\ref{C_i*_D_i*_bound}))}  \\
&> & (n/10)^{k-1} \mbox{ (since $\sqrt{\varepsilon} < 1/(100k^2)$, $\tau = 1/(9k)$ and $n\ge 100k^2$)}\\
&> & \sqrt{\varepsilon} n^{k-1} \mbox{ (since $\sqrt{\varepsilon} < 1/(k(10k^2)^{k-1})$},
\end{eqnarray*}
contradicting the fact that $v \not\in N$.

Let $M^j \subseteq T^j$ with $|M^j| = \lfloor \tau n\rfloor$. Thus, this process works  for all $j\in [k]$, and we see that
$M^0,M^1,\ldots, M^k$ give the desired matchings for Claim 6. $\Box$

\medskip

By Claim 6, there exist pairwise disjoint matchings $M^0, M^1, \ldots, M^k$ in $H^*$
such that
\begin{itemize}
\item $M^i$ is of $i$-type for $i\in [k]\cup \{0\}$, and
\item $|M^i| = |M^1|\ge \lfloor \tau n\rfloor$ for $i\in [k]$.
\end{itemize}
We choose such $M^0, M^1, \ldots, M^k$ that
\begin{itemize}
\item $|M^1|=|M^2|=\ldots =|M^k|$ is maximum and, subject to this,
\item  $|M^0|$ is maximum.
\end{itemize}
Let $M=\bigcup_{i\in [k]\cup \{0\}} M^i$. By Claim 6, we have,  for $i \in [k]$,
$$|D_i^* \cap V(M)| \equiv 0 \pmod {k-1} \mbox{ and } |M^i|\le |D_i^*|/(k-1).$$
\medskip

\textit{Claim 7.}  $|M^0| \geq \tau n$.

For, otherwise, suppose $|M^0| < \tau n$.
Note that  for $i \in [k]$,
\begin{eqnarray*}
& & |C_i^* - V(M^0 \cup M^1 \cup \cdots \cup M^k)|  \\
&=& |C_i^*| - |M^0| - |M^i| \\
&> & |C_i^*| - \tau n - |D_i^*|/(k-1) \quad \mbox{ (since $|M^0| < \tau n$ and $|M^i|\le |D_i^*|/(k-1)$)}\\
&\geq & |C_i^*| - \tau n - (n -|C^*_i|)/(k-1)  \\
&= & k|C_i^*|/(k-1) - \tau n - n/(k-1) \\
&\geq & ((1/2 - 2\sqrt{\varepsilon}k)n - 2 - 8 \sqrt{\varepsilon} k^2 n - k^2/2)k/(k-1) - \tau n - n/(k-1) \quad\mbox{ (by (\ref{C_i*_D_i*_bound})) } \\
&> & n/10  \mbox{ (since $\sqrt{\varepsilon} < 1/(100k^2)$, $\tau  = 1/(9k)$ and $n\ge 100k^2$).}
\end{eqnarray*}

Thus there exists $v \in C_1^*  - V(M^0 \cup M^1 \cup \cdots \cup M^k)$.
Since $|M^0|$ is maximized,
$C^*- V(M^0 \cup M^1 \cup \cdots \cup M^k)$ has no edge. Therefore,
\begin{eqnarray*}
|N_{H_0(k,n)}(v) - N_H(v)|
&\geq &\prod_{i \in [k] - \{1\}} |C_i^* - V(M^0 \cup M^1 \cup \cdots \cup M^k)| \\
&> &( n/10 )^{k-1}\\
&> &\sqrt{\varepsilon} n^{k-1} \mbox{ (since $\sqrt{\varepsilon} < 1/(k(10k^2)^{k-1})$},
\end{eqnarray*}
contradicting the fact that $v \not\in N$. $\Box$

\medskip

\textit{Claim 8.}  $D^* \subseteq V(M)$.

For, otherwise, suppose that $D^* - V(M) \neq \emptyset$.
Recall that for each $j$-type edge $f$, $|f \cap C^*| = |f \cap C_j| = 1$.
Since $|D_i^*| \equiv 0 \pmod {k-1}$ (by Claim 5) and
$|D_i^* \cap V(M)| \equiv 0 \pmod {k-1}$, it follows that
$|D_i^* - V(M)|\ge k-1$ for $i \in [k]$.
So, for $i \in [k]$,
let $s_{i,1},s_{i,2},\ldots,s_{i,k-1}\in D_i^* - V(M)$ be
distinct.

When $C_i^*-V(M)\ne \emptyset$ for $i\in [k]$,
let $w_i\in C_i^* - V(M)$ for $i \in [k]$;  otherwise let $\{w_1,\ldots,w_k\} \in M^0$ with $w_i \in C_i^*$ for $i \in [k]$ (by Claim 7).
Let $S_j :=\{w_j, s_{i,j}: i\in [k]-\{j\}\}$ for $j\in [k-1]$, and let
$S_k := \{w_k, s_{i,i}: i\in [k-1]\}$.

Suppose for each $j \in [k]$ there exist distinct
$e_1^j,\ldots,e_{k-1}^j\in M^j$ such that $H^*[e_1^j \cup \cdots \cup
e_{k-1}^j \cup S_j]$ contains a perfect matching $\{f_1^j, \ldots,
f_k^j\}$. Then, $N^j:=(M^j-\{e_1^j,\ldots,e_{k-1}^j\})\cup
\{f_1^j,\ldots,f_k^j\}$ is a  matching in $H^*$ for each $j \in [k]$, and
$|N^j|=|N^1|>\lfloor \tau n\rfloor$ for $j\in [k]$. Let $N^0=M^0-\{\{w_1,\ldots, w_k\}\}$. Then
 $N^0,N^1,\ldots, N^k$ are pairwise disjoint. However, $|N^j|=|M^j|+1$
 for $j\in [k]$, contradicting the choice of $M^0, M^1,\ldots, M^k$.

Thus we may assume without loss of generality that
for any $k-1$ distinct edges $e_1^k,\ldots,e_{k-1}^k \in M^k$,  $H^*[e_1^k \cup \cdots \cup
e_{k-1}^k \cup S_k]$ has no perfect matching.
For $i \in [k-1]$, let $e_i^k := \{v_{i,1},v_{i,2},\ldots,v_{i,k}\}$ with $v_{i,k} \in C_k^*$ and $v_{i,j} \in D_j^*$ for $j \in [k-1]$.
For convenience, let  $v_{k,k} := w_k$ and $v_{k,j} := s_{j,j}$ for $j \in [k-1]$.
For $i \in [k]$, define $f_i^k := \{v_{1,i+1},
v_{2,i+2},\ldots,v_{k-1,i+k-1},v_{k,k+i}\}$, where the addition in the
subscripts is modulo $k$ (except that we write $k$ for $0$).
Then $f_i^k \not\in E(H^*)$ for some $i \in [k]$, as otherwise, $\{f_1^k, \ldots, f_k^k\}$ would be a perfect matching in $H^*
[e_1^k \cup \cdots \cup e_{k-1}^k \cup S_k]$.
 Since  $e_1^k,\ldots,e_{k-1}^k\in M^k$ are chosen arbitrarily and $k$ is odd (by Claim 5), we have
\begin{eqnarray*}
& &  |\{e \in E(H_0(k,n)) - E(H): |e\cap  \{v_{k,i}:i\in [k]\}|=1\}|\\
&\geq & \binom{|M^k|}{k-1}\\
&= &|M^k||M^k - 1|\cdots |M^k - (k-1) + 1|/(k-1)! \\
&> & {((\lfloor \tau n\rfloor - k)/(k-1))}^{k-1} \\
&> & (n/(10k^2))^{k-1} \mbox{ (since $\tau  = 1/(9k)$ and $n\ge 100k^2$)}  \\
&> & k\sqrt{\varepsilon} n^{k-1}  \mbox{ (since $\sqrt{\varepsilon} < 1/(k(10k^2)^{k-1})$}.
\end{eqnarray*}
So there exists $i \in [k]$ such that
$|N_{H_0(k,n)}(v_{k,i}) - N_H(v_{k,i})| > \sqrt{\varepsilon} n^{k-1}$,
contradicting the fact that $v_{k,i} \not\in N$. $\Box$

\medskip

If $C^*\subseteq V(M)$ then, by Claim 8, $M$ is a perfect matching in $H^*$; so
$\{e_0\} \cup M_1 \cup M_2 \cup M_3 \cup M$ is a perfect matching in
$H$.

Therefore, we may assume that  $C^*\not\subseteq V(M)$, and let $w_i \in C_i^* - V(M)$ for $i \in [k]$.
Note that $|M^0| \ge \tau n > k - 1$  (by Claim 7).
Let $e_1,\ldots,e_{k-1} \in M^0$ be distinct and chosen arbitrarily.
Let $e_i := \{v_{i,1},v_{i,2},\ldots,v_{i,k}\}$ for $i \in [k-1]$,
where $v_{i,j} \in C_j^*$ for $j \in [k]$.
For $i \in [k]$, define $f_i := \{w_i, v_{1,i+1},
v_{2,i+2},\ldots,v_{k-1,i+k-1}\}$, with the addition in the subscripts
taken modulo $k$ (except we use $k$ for $0$).

If $f_i \in E(H^*)$ for all $i \in [k]$, then $N^0:= (M^0 \cup
\{f_1,\ldots,f_k\}) - \{e_1,\ldots,e_{k-1}\}$ is a matching
in $H^*$ with $|N^0| = |M^0| + 1$; so $N^0, M^1, \ldots, M^k$
contradict the choice of $M^0,M^1,\ldots, M^k$.

Hence, $f_i \not\in E(H^*)$ for some $i \in [k]$.
Since $e_1,\ldots,e_{k-1}\in M^0$ are chosen arbitrarily and $k$ is odd, we have
\begin{eqnarray*}
& & |\{e \in E(H_0(k,n)) - E(H): |e\cap \{w_i: i\in [k]\}|=1\}| \\
&\geq & \binom{|M^0|}{k-1}\\
&= & |M^0||M^0 - 1|\cdots |M^0 - (k-1) + 1|/(k-1)! \\
&>& {((\lfloor \tau n\rfloor - k)/k)}^{k-1} \\
&> &(n/(10k^2))^{k-1}  \mbox{ (since $\tau  = 1/(9k)$ and $n\ge 100k^2$)}\\
&> & k \sqrt{\varepsilon} n^{k-1}  \mbox{ (since $\sqrt{\varepsilon} < 1/(k(10k^2)^{k-1})$}.
\end{eqnarray*}
So there exists $i \in [k]$ such that
$|N_{H_0(k,n)}(w_i) - N_H(w_i)| > \sqrt{\varepsilon} n^{k-1}$,
contradicting the fact that $w_i \not\in N$. \qed

\begin{coro}
\label{near_critical3}
Let $k \geq 3$ be a positive integer, and let   $\varepsilon >0$ be such
 that $\sqrt{\varepsilon}< \min\{1/(100  k^2),$ $1/(k (10k^2)^{k-1})
\}$.
Let $H$ be a $k$-partite $k$-graph with $n > 100k^2$ vertices in each
partition class, such that  $\delta_{k-1}(H) \geq \lfloor n/2\rfloor$
and  $H$ is $\varepsilon$-close to $H_0(k,n)$.
Then $H$ has no perfect matching if, and only if,
\begin{itemize}
\item[$(i)$] $k$ is odd, $n\equiv 2\pmod 4$, and $H\cong H_0(k,n)$, or
\item [$(ii)$] $n$ is odd and there exist $d_i\in \{(n+1)/2,(n-1)/2\}$
  for $i\in [k]$ such that  $\sum_{i=1}^k d_i$ is odd and $H \subseteq H_0(d_1,\ldots,d_k;k,n)$.
\end{itemize}
\end{coro}
\pf
Let $H$ be a $k$-partite $k$-graph with $n$ vertices in each partition
class, such that
 $\delta_{k-1}(H) \geq \lfloor n/2\rfloor$
and $H$ is $\varepsilon$-close to $H_0(k,n)$.

Suppose  $(i)$ or $(ii)$ holds. Then there exist integers $d_1,\ldots,
d_k$ such that $\sum_{i=1}^k d_i$ is odd,
$H\subseteq H_0(d_1,\ldots,d_k;k,n)$,  $d_1 = \ldots = d_k = n/2$ when
$(i)$ holds, and  $d_i \in \{(n+1)/2, (n-1)/2\}$ when $(ii)$ holds.
By the definition of $H_0(d_1,\ldots,d_k;k,n)$,
there exists $D \subseteq V(H)$ such that
$|D| = \sum_{i=1}^k d_i$  is odd and
$|e\cap D|$ is even for all $e\in E(H)$.
Hence, $H$ contains no perfect matching.

 Next, suppose $H$ has no perfect matching.
 Applying Lemma \ref{near_critical} with $\alpha = 1/8$, we may assume
 that there exist  $d_i\in [ \lceil 3n/8 \rceil, \lfloor 5n/8 \rfloor
 ]$ for $i \in [k]$ such that $\sum_{i=1}^k d_i$ is odd and $H
 \subseteq H_0(d_1,\ldots,d_k;k,n)$. Let $V_1,\ldots, V_k$ be the
 partition classes of $H$ and $H_0(d_1,\ldots,d_k;k,n)$.
 For $i \in [k]$, let $D_i \subseteq V_i$ be such that $|D_i| = d_i$ and $|e \cap (\cup_{j \in [k]} D_j)|$ is even for all $e \in E(H)$.

We claim that
$\delta_{k-1}(H)\leq \min \{d_i,n-d_i\}$ for all $i \in [k]$. By
symmetry, we only show $\delta_{k-1}(H)\le \min\{d_1, n-d_1\}$.
Let $S :=\{v_2,\ldots,v_k\}$ be a legal set such that $v_2 \in D_2$
and $v_i\in V_i - D_i $ for $i \in [k]- \{1,2\}$;
then, since $e \cap D_1 \neq \emptyset$ for all $e \in E(H)$ with $S \subseteq e$,
$\delta_{k-1}(H)\leq d_{H}(S)\leq |D_1| = d_1$.
Let $T :=\{u_2,\ldots,u_k\}$ be a legal set such that $u_i\in V_i -
D_i $ for $i \in [k]- \{1\}$; then, since $e\cap D_1=\emptyset$ for
any $e\in E(H)$ with $T\subseteq e$,
$\delta_{k-1}(H)\leq d_{H}(T)\leq |V_1 - D_1| = n-d_1$.
Hence,  $\delta_{k-1}(H)\leq \min\{d_1,n-d_1\}$.

If $n$ is odd then $\delta_{k-1}(H)\geq\lfloor n/2\rfloor=(n-1)/2$;
so  by the above claim, $d_i\in \{(n-1)/2,(n+1)/2\}$ for all $i\in [k]$, and $(ii)$ holds.
Thus, we may assume that $n$ is even.
Then by the above claim, $d_i=n/2$ for all $i\in [k]$.
Recall that $\sum_{i=1}^k d_i$ is odd.
Thus both $n/2$ and $k$ are odd, and hence $n\equiv 2 \pmod 4$.
Since $H\subseteq H_0(d_1,\ldots,d_k;k,n) =  H_0(k,n)$, we have
$H=H_0(k,n)$ and $(i)$
holds.
\qed

\section{Hypergraphs not close to $H_0(k,n)$}

In this section, we prove Theorem~\ref{main} for  hypergraphs that are not close  to $H_0(k,n)$, see Lemma~\ref{lemma_away_from_extremal}. For this, we need a result on almost perfect
matchings in $k$-partite $k$-graphs.

K\"uhn and Osthus \cite{Kuh06}  showed that if $H$ is a $k$-partite
$k$-graph  with each partition
classes of size $n$ and $\delta_{k-1}(H)\geq n/k$, then  $H$ has a matching of size at least $n-(k-2)$.
R\"odl and Ruci\'nski \cite{Rod09} asked the following question: Is it
true that
 $\delta_{k-1}(H) \geq n/k$ implies that $H$ has a matching of size at least $n-1$?
The present authors \cite{Lu16} and,
independently, Han, Zang, and Zhao \cite{HZZ16} answered this question affirmatively for large
$n$.
\begin{lemma}
\label{near_perfect}
Let $k,n$ be positive integers with $k \geq 3$ and $n$ sufficiently
large,  and let $H$ be a $k$-partite $k$-graph with $n$ vertices in each partition class.
If $\delta_{k-1}(H) \geq n/k$, then $H$ has a matching of size at least $n-1$.
\end{lemma}

Let $k \geq 2$ be a positive integer
and $H$ be a  $k$-partite $k$-graph with partition classes $V_1,\ldots,V_k$.
Given $N_i\subseteq V_i$  for $i\in [k]$, let
$$E_H(N_1,\ldots,N_k):=\{ e\in E(H) : e\subseteq \cup_{i\in [k]}
N_i \},$$
and
$$e_H(N_1,\ldots,N_k):=|E_H(N_1,\ldots,N_k)|.$$
For $j\in[k]$, let $$\Lambda_j := \{ \{v_i\in V_i : i\in [k]-\{j\}\}:  d_H(\{v_i:i\in [k]-\{j\}\}) \geq (1/2+2/\log n)n \}.$$

\begin{lemma}\label{Large-Set}
Let $k \geq 2$ be a positive integer.
For any $\varepsilon>0$, there exists $n_0 > 0$ such that the following holds.
Let $H$ be a $k$-partite $k$-graph with partition classes $V_1,\ldots,V_k$ such that $|V_i| = n\geq n_0$ for $i \in [k]$
and $\delta_{k-1}(H)\geq (1/2-1/\log n )n$.
Suppose $H$ is not $\varepsilon$-close to $H_0(k,n)$.
Then one of the following conclusions holds:

\begin{itemize}
\item [$(i)$]  For all  $i\in [k]$ and $N_i\subseteq V_i$ with $|N_i|\geq (1/2-1/\log
  n)n$, $e_H(N_1,\ldots,N_k)\geq n^k/\log^3n$.

\item [$(ii)$] There exists $j\in [k]$ such that $|\Lambda_j|\geq
  n^{k-1}/\log n$.
\end{itemize}
\end{lemma}

\pf
Let $H$ be a $k$-partite $k$-graph with partition classes $V_1,\ldots,V_k$ such that $|V_i| = n$ for $i \in [k]$.
For convenience, let $\gamma:=1/\log n$. Then $\delta_{k-1}(H)\geq (1/2- \gamma)n$.

Suppose $H$ is not $\varepsilon$-close to $H_0(k,n)$, and assume that  neither $(i)$ nor $(ii)$  holds.
Then there exist $N_1,\ldots,N_k$ with $N_i\subseteq V_i$ and
$|N_i|\geq (1/2-\gamma)n$ for $i\in [k]$ such that
\begin{align}\label{N1-Nk-up}
e_H(N_1,\ldots,N_k)<\frac{n^k}{\log^3n} = o(n^k),
\end{align}
and, for all $j \in [k]$,
\begin{align}\label{Lambda-low}
|\Lambda_j|< \frac{n^{k-1}}{\log n} = \gamma n^{k-1}.
\end{align}

\medskip

{\it Claim~1.}  $|N_i|<(1/2+2\gamma)n$ for  $i \in [k]$.

For, otherwise, we may assume without loss of generality  that
$|N_{k}|\geq (1/2+2\gamma)n$. Then $|V_k-N_k|\le (1/2-2\gamma)n$.
 For any legal $(k-1)$-set
$\{v_1,\ldots, v_{k-1}\}$ with $v_i\in N_i$ for $i \in [k-1]$,  we
have
\[
|N_H(v_1,\ldots,v_{k-1})\cap N_k|\geq \delta_{k-1}(H)-|V_k-N_k| \ge (1/2-\gamma)n-  (1/2-2\gamma)n = \gamma n.
\]
Hence, by choosing $n_0$ large enough, we have for $n\ge n_0$,
\[
e_H(N_1,\ldots,N_k)\geq |N_1|\cdots |N_{k-1}||N_H(\{v_1,\ldots, v_{k-1}\})\cap N_k|\ge ((1/2-\gamma)n)^{k-1} \gamma n >\frac{n^k}{\log^3 n},
\]
contradicting (\ref{N1-Nk-up}). $\Box$

\medskip

For $i\in [k]$,
let $N_i':=V_i-N_i$ and $A_i \in \{N_i,N_i'\}$. Since
$|N_i|\ge (1/2-\gamma)n$, $|N_i'|\le (1/2+\gamma)n$. By Claim 1,
$|N_i'|> (1/2-2\gamma)n$. Therefore, for $i\in [k]$,
\begin{align}\label{Ai-upper-lower-bound}
(1/2 - 2\gamma)n < |A_i| \leq (1/2 + 2\gamma)n.
\end{align}

\medskip

{\it Claim~2.} For $i\in [k]$,
$e_H(A_1,\ldots,A_{i-1}, V_i,A_{i+1},\ldots, A_k)=(n/2)^k+o(n^k)$.

By symmetry, we only  prove Claim 2 for the case when  $i=k$.
Note that
\begin{align*}
e_H(A_1,\ldots,A_{k-1},V_k)
&\geq \left( \prod_{i=1}^{k-1}|A_i|\right)  (1/2-\gamma)n   \quad \mbox{ (since $\delta_{k-1}(H) \geq (1/2 - \gamma)n$)} \\
&\geq \left((1/2-2\gamma)n\right)^{k-1} (1/2-\gamma)n  \quad \mbox{  (by (\ref{Ai-upper-lower-bound})) } \\
&=(n/2)^k+o(n^k).
\end{align*}
On the other hand,
\begin{align*}
e_H(A_1,\ldots,A_{k-1}, V_k)
&\leq |\Lambda_k| n+\left(\prod_{i=1}^{k-1}|A_j|\right)(1/2+2\gamma)n\\
&< \gamma n^k +\left((1/2+2\gamma)n\right)^{k-1}(1/2+2\gamma)n  \quad
  \mbox{(by    (\ref{Lambda-low}) and  (\ref{Ai-upper-lower-bound}))}\\
&=(n/2)^k+o(n^k). \quad \Box
\end{align*}

\medskip

{\it Claim 3}.   Let $I(A_1,\ldots, A_k) :=\{i\in [k]\ : \
A_i=N_i'\}$.
\begin{itemize}
\item [($i)$] If $|I(A_1,\ldots, A_k)|$ is odd
then $e_H(A_1,\ldots,A_k)=(n/2)^k+o(n^k)$, and
\item [($ii$)]  if $|I(A_1,\ldots,A_k)|$ is even then $e_H(A_1,\ldots,A_k)=o(n^k)$.
\end{itemize}

We apply induction on $|I(A_1,\ldots, A_k)|$. When $|I(A_1,\ldots, A_k)|=0$, we have
$e_H(A_1, \ldots,$ $A_k) = e_H(N_1,\ldots ,N_k) = o(n^k)$ by  (\ref{N1-Nk-up}). When
$|I(A_1,\ldots, A_k)|=1$,
say $A_i = N_i'$ and $A_j = N_j$ for $j \in [k] - \{i\}$,
then we have
\begin{eqnarray*}
& &e_H(A_1,\ldots,A_k)\\
& = & e_H(A_1,\ldots, A_{i-1}, V_i, A_{i+1}, \dots, A_k) -
  e_H(N_1,\ldots, N_{i-1}, N_i, N_{i+1}, \dots, N_k)\\
& = & (n/2)^k+o(n^k)
\end{eqnarray*}
 by Claim 2 and
(\ref{N1-Nk-up}).

Now assume Claim 3 holds for $A_1,\ldots, A_k$ with $A_i\in \{N_i,N_i'\}$ and $0\le
|I(A_1,\ldots, A_k)|=l<k$.  Consider a choice of  $A_i\in \{N_i,N_i'\}$ for $i\in
[k]$ with $|I(A_1,\ldots, A_k)|=l+1$.
Let $A_j = N_j'$ for some $j \in [k]$. Observe that
$$e_H(A_1,\ldots,A_k) = e_H(A_1,\ldots, A_{j-1}, V_j, A_{j+1}, \dots, A_k) - e_H(A_1,\ldots, A_{j-1}, N_j, A_{i+1}, \dots, A_k).$$
Therefore, by  (\ref{N1-Nk-up}) and Claim 2, it follows from the induction hypothesis that
if $l+1$ is odd then $l$ is even and  $e_H(A_1,\ldots,A_k) =
((n/2)^k+o(n^k)) - o(n^k)=(n/2)^k+o(n^k)$, and
if $l+1$ is even then $l$ is odd and
 $e_H(A_1,\ldots,A_k) =((n/2)^k+o(n^k)) - ((n/2)^k+o(n^k))=o(n^k)$. $\Box$

\medskip

For $i\in [k]$,
let $B_i\subseteq V_i$ be such that $|B_i|=\lfloor n/2\rfloor$ with
$|B_i\cap N_i|$  maximal,
and let $B_i':=V_i-B_i$.
By (\ref{Ai-upper-lower-bound}), for $i\in [k]$,
$|B_i-N_i|\leq 2\gamma n$ and  $|B_i'-N_i'|\leq 2\gamma n$. Hence, if
$A_i\in \{N_i,N_i'\}$ and $C_i\in \{B_i,B_i'\}$ such that  for $i\in
[k]$,
$A_i=N_i$  iff $C_i=B_i$, then
$$\left|E_{H}(A_1,\ldots,A_k)-E_{H}(C_1,\ldots,C_k)\right|=o(n^k)$$
and
\begin{eqnarray*}
& & \left|E_{H_0(k,n)}(C_1,\ldots,C_k)-E_{H}(C_1,\ldots,C_k)\right|\\
&\le & \left|E_{H_0(k,n)}(C_1,\ldots,C_k)-E_{H}(A_1,\ldots,A_k)\right|+ \left|E_H(A_1,\ldots,A_k)-E_{H}(C_1,\ldots,C_k)\right|\\
&=& \left|E_{H_0(k,n)}(C_1,\ldots,C_k)-E_{H}(A_1,\ldots,A_k)\right|+o(n^k).
\end{eqnarray*}

For $i \in [k]$, let $B_i$ play the role of $D_i$ in the
definition of $H_0(k,n)$.
Then, for any $\varepsilon > 0$,
\begin{eqnarray*}
& & |E(H_0(k,n))-E(H)|\\
&\leq & \sum_{C_i \in \{B_i,B_i'\}, i \in [k] } \left|E_{H_0(k,n)}(C_1,\ldots,C_k)-E_{H}(C_1,\ldots,C_k)\right|\\
&\leq & \sum_{C_i \in \{B_i,B_i'\}, A_i\in \{N_i,N_i'\} \atop A_i =
        N_i \text{ iff } C_i = B_i \text{ for } i \in [k] } \left(\left|E_{H_0(k,n)}(C_1,\ldots,C_k)-E_{H}(A_1,\ldots,A_k)\right| + o(n^k) \right)\\
&\leq & \sum_{C_i\in \{B_i,B_i'\} \text{ for } i\in [k]} \left( \left( \sum_{i \in [k]} \sum_{v \in (B_i - N_i) \cup (B_i' - N_i')} |N_{H_0(k,n)}(v)|\right) + o(n^k) \right) \\
&\leq & 2^{k} \left( k (4\gamma n)  n^{k-1} + o(n^k) \right) \quad
        \mbox{ (since $|B_i - N_i| \le 2\gamma n$ and $|B_i' - N_i'|
        \le 2\gamma n$)}\\
&\leq &\varepsilon n^k \quad \mbox{ (since $\gamma=1/\log n$ and we
        may choose $n_0$ large enough)} .
\end{eqnarray*}
However, this contradicts the assumption that $H$ is not $\varepsilon$-close to $H_0(k,n)$. \qed

\medskip

Next, we define two ``absorbing'' matchings for a legal $k$-set $S$
in a $k$-partite $k$-graph. This concept was  first considered by R\"odl, Ruci\'nski, and Szemer\'edi \cite{RRS09}.
Let $k \geq 3$ be a positive integer and $H$ be a $k$-partite
$k$-graph.

\begin{figure}[!htb]
\centering
\includegraphics[width=8cm]{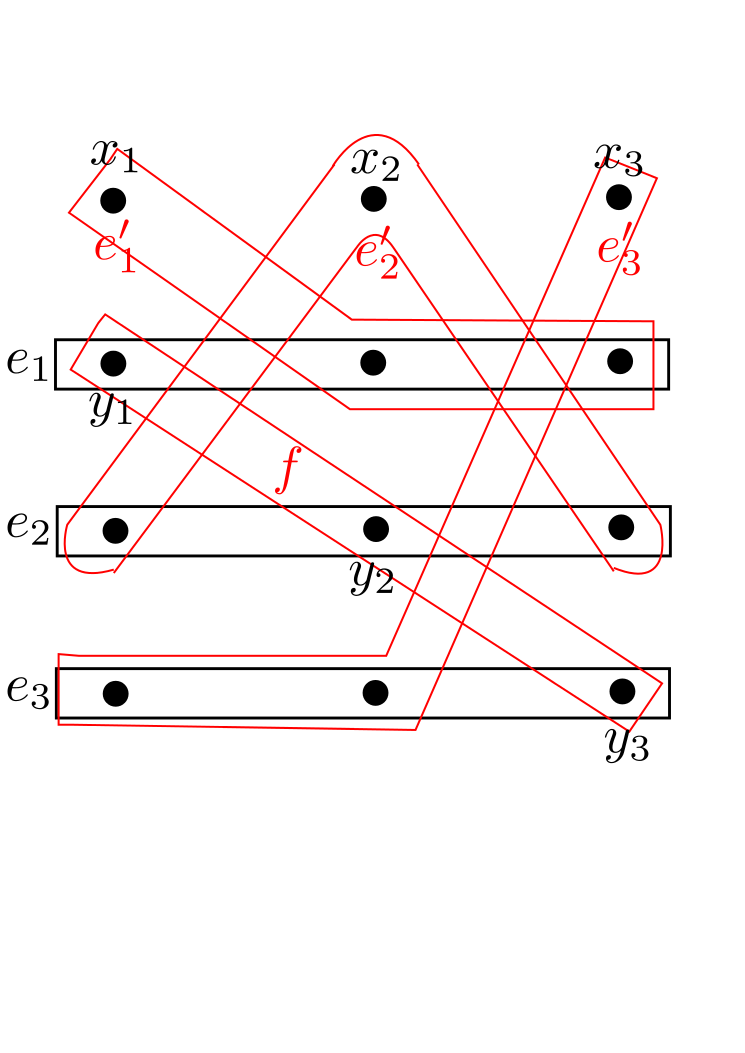}
\caption{$\{x_1,x_2,x_3\}$-absorbing matchings}
\label{fig:gull}
\end{figure}

Given a legal $k$-set $S=\{x_1,\ldots, x_k\}$ in a $k$-partite $k$-graph
$H$, a $k$-matching $\{e_1,\ldots, e_k\}$ in $H$ is said to be
\textbf{$S$-absorbing} if
there is  a $(k+1)$-matching $\{e_1',\ldots,e_k',f\}$ in $H$ with
$f=\{y_1,\ldots,y_k\}$ such that
\begin{itemize}
\item $e_i'\cap e_j=\emptyset$ for all $i\neq j$,
\item $e_i'-e_i=\{x_i\}$ and $e_i-e_i'=\{y_i\}$ for  $i\in [k]$.
\end{itemize}
Figure \ref{fig:gull} illustrates  an $\{x_1,x_2,x_3\}$-absorbing 3-matching $\{e_1,e_2,e_3\}$.

\begin{figure}[!htb]
\centering
\includegraphics[width=8cm]{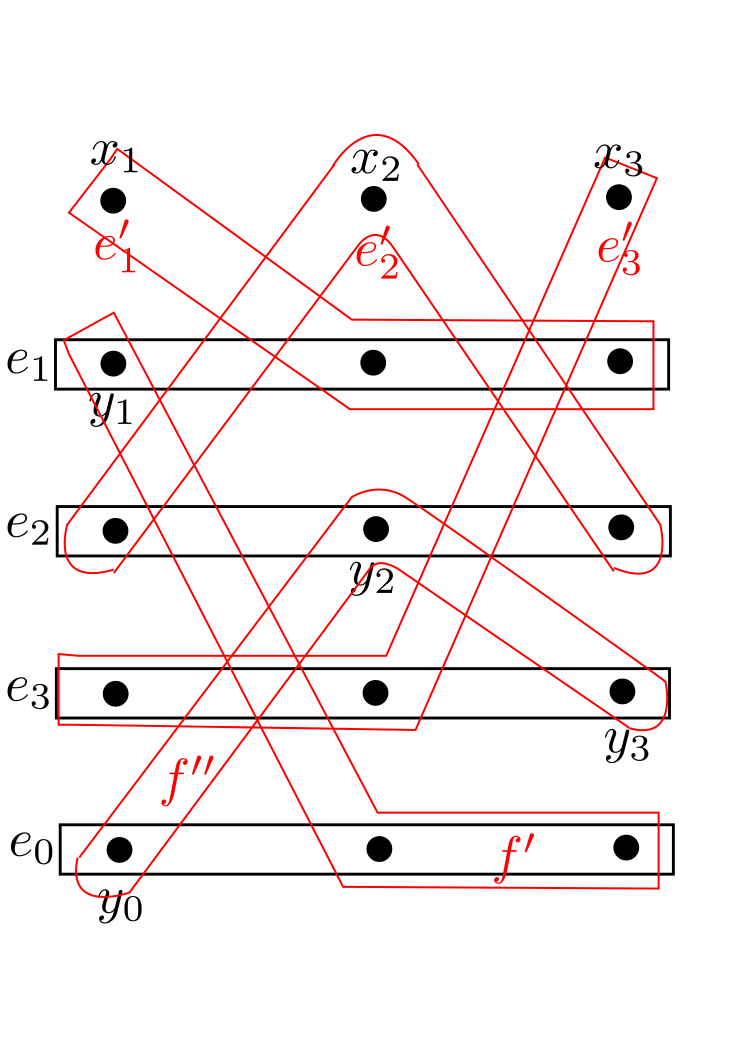}
\caption{$\{x_1,x_2,x_3\}$-absorbing $(k+1)$-matching  for $k=3$}
\label{fig:gull2}
\end{figure}

Given a legal $k$-set $S=\{x_1,\ldots, x_k\}$ in a $k$-partite
$k$-graph $H$,  a $(k+1)$-matching $\{e_0,e_1,\ldots, e_k\}$ in $H$ is
said to be \textbf{$S$-absorbing} if there is  a $(k+2)$-matching
$\{e_1',\ldots,e_k',f',f''\}$ in $H$, with $e_1\cap
f'=f'-e_0=\{y_1\}$, $e_0-f'=\{y_0\}$, and
$f'':=\{y_0,y_2,\ldots, y_k\}$,  such that
\begin{itemize}
\item $e_i'\cap e_j=\emptyset$ for all $i\neq j$, and
\item $e_i'-e_i=\{x_i\}$ and $e_i-e_i'=\{y_i\}$ for all $i\in [k]$.
\end{itemize}
Figure  \ref{fig:gull2} illustrates an $\{x_1,x_2,x_3\}$-absorbing 4-matching $\{e_0,
e_1,e_2,e_3\}$.

\medskip

The next result says that no matter which conclusion of Lemma~\ref{Large-Set}
holds, there are always many $S$-absorbing matchings in $H$ for any
given legal $k$-set $S$.

\begin{lemma}\label{absorb-num}
Let $k \geq 3$ be a positive integer.
There exists $n_1 > 0$ such that the following holds. Let $H$ be a
$k$-partite $k$-graph with $n\geq n_1$ vertices in each partition
class and with $\delta_{k-1}(H)\geq (1/2-1/\log n)n$.
Let $S:=\{x_1,\ldots,x_k\} \subseteq V(H)$ be legal.
\begin{itemize}

\item [$(i)$]  If  for all  $i\in
  [k]$ and $N_i\subseteq V_i$ with  $|N_i|\geq (1/2-1/\log n)n$,  we have
  $e_H(N_1,\ldots,N_k)\geq n^k/\log^3n$,
then the number of $S$-absorbing $k$-matchings in $H$ is $\Omega(n^{k^2}/\log^3 n)$.

\item [$(ii)$] If there exists $j\in [k]$ such that $|\Lambda_j|\geq n^{k-1}/\log n$,
then the number of $S$-absorbing $(k+1)$-matchings in $H$ is $\Omega(n^{k^2+k}/\log^3 n)$.
\end{itemize}
\end{lemma}

\pf To prove $(i)$, we assume that,  for all $i\in [k]$ and
  $N_i\subseteq V_i$ with $|N_i|\geq (1/2-1/\log
  n)n$, we have  $e_H(N_1,\ldots,N_k)\geq
n^k/\log^3n$.

  Note that, for each $i\in [k]$, $x_i$ is contained in
$(n-1)^{k-2}$ legal $(k-1)$-sets in $H$ that are disjoint from $S$ and
one given partition class of $H$,
and each such legal $(k-1)$-set is contained in at least  $(1/2-1/\log
n)n-1$ edges in $H-S$ (since $\delta_{k-1}(H)\geq (1/2-1/\log n)n$).
Thus, there exists $n_1$ such that if $n \geq n_1$, there are at least $ n^{k-1}/3$ legal $(k-1)$-sets $B_i$ disjoint from $S$ such  that
$e_i':=\{x_i\}\cup B_i\in E(H)$.

By a similar argument (and choosing $n_1$ large enough),
there are at least $((n-k)^{(k-1)}/3)^k\ge (1/3-o(1))^kn^{(k-1)k}$ choices
of pairwise disjoint such legal $(k-1)$-sets $B_1,\ldots,B_k$.

For each such choice of $B_1,\ldots, B_k$, let $N_i :=N_H(B_i)$ for $i\in[k]$.
 Then $|N_i|\geq \delta_{k-1}(H) \geq
 (1/2-1/\log n)n$.  By assumption, $e_H(N_1,\ldots, N_k)\ge
 n^k/\log^3n$; so there are at least $n^k/\log^3 n - k^2 n^{k-1}$ choices
 of an edge $f :=\{y_1,\ldots,y_k\}$ from $H[\cup_{i \in [k]} N_i] -
 \cup_{i \in [k]} e_i'$ such that
$e_i:=B_i\cup \{y_i\}\in E(H)$ for $i \in [k]$.

Hence, the number of $S$-absorbing $k$-matchings $\{e_1,\ldots,e_k\}$ is at least
\[
(1/3 - o(1))^kn^{(k-1)k} (n^k/\log^3n - k^2 n^{k-1}) = \Omega(n^{k^2}/\log^3 n),
\]
as claimed in $(i)$.

\medskip

We now prove $(ii)$. So assume without loss of generality that  $|\Lambda_1|\geq
  n^{k-1}/\log n$.
As in the previous case, since $\delta_{k-1}(H)\geq (1/2-1/\log n)n$,
there are at least $(1/3-o(1))^kn^{(k-1)k} = \Omega(n^{(k-1)k})$
choices of disjoint legal $(k-1)$-sets $B_1,\ldots,B_k$ such that
$\{x_i\}\cup B_i\in E(H)$ for $i\in [k]$.

For $i=2,\ldots,k$, we choose $y_i\in N_H(B_i) - \{x_i\}$ and let
$e_i:=B_i\cup \{y_i\}$. Note that we have  $(1/2-2/\log n)n - 1 =
\Omega(n)$ choices for each $y_i$.

By assumption, there are at least $n^{k-1}/\log n - k(k+2)n^{k-2} =
\Omega(n^{k-1}/\log n)$ choices for a $(k-1)$-set $T\in \Lambda_1$
that is disjoint from
$S\cup B_1\cup\cdots\cup B_k\cup\{y_2,\ldots,y_k\}$.
Since $N_H(T)>(1/2+2/\log n)n$ (as $T \in \Lambda_1$) and
$\delta_{k-1}(H) \geq (1/2 - 1/\log n)n$, we have $|N_H(T)\cap N_H(B_1)|\geq n/\log n$ and $|N_H(T)\cap N_H(\{y_2,\ldots,y_k\})|\geq n/\log n$.
Consequently, there exist
distinct $y_0$ and $y_1$ with $y_1\in (N_H(T)\cap N_H(B_1))-(S\cup B_1\cup\cdots\cup B_k)$ and $y_0\in
(N_H(T)\cap N_H(\{y_2,\ldots,y_k\}))-(S\cup B_1\cup\ldots\cup B_k)$,
and  there are at least $n/\log n-k(k+1) - 1$ choices for each of $y_0$ and
$y_1$.

Let $e_0:=\{y_0\} \cup T$, $e_1 := \{y_1\} \cup B_1$,
$f':=\{y_1\} \cup T$ and  $f'':= \{y_0,y_2,y_3,\ldots,y_k\}$. Then  $\{e_0,\ldots,e_k\}$ is an $S$-absorbing $(k+1)$-matching (using $e_i'=B_i\cup \{x_i\}$ for
$i\in [k]$).
Moreover, the number of choice for $\{e_0,\ldots,e_k\}$ is the product of the numbers of choices for $B_1,\ldots,B_k$,
$y_2,\ldots,y_k$, $T$, $y_0$, $y_1$, which is at least
\[
\Omega\left(n^{k(k-1)}\right)  \Omega \left(n^{k-1} \right)
\Omega \left(\frac{n^{k-1}}{\log n} \right) \left(\frac{n}{\log n}-k(k+1)-1 \right)^2
= \Omega \left(\frac{n^{k^2+k}}{\log^3 n} \right).
\]
So we have $(ii)$. \qed

\medskip

We will need to use  Chernoff bounds, which can be found in \cite{mitzenmacher}. 

\begin{lemma}
\label{chernoff}
Suppose $X_1, ..., X_n$ are independent random variables taking values in $\{0, 1\}$. Let $X$ denote their sum and $\mu = \mathbb{E}[X]$ denote the expected value of $X$. Then for any $0 < \delta \leq 1$,
$$\mathbb{P}[X \geq (1+\delta) \mu] < e^{-\frac{\delta^2
    \mu}{3}}\mbox{ and } \mathbb{P}[X \leq (1-\delta) \mu] < e^{-\frac{\delta^2 \mu}{2}},$$
and for any $\delta \geq 1$,
$$\mathbb{P}[X \geq (1+\delta) \mu] < e^{-\frac{\delta \mu}{3}}.$$
\end{lemma}

We now show that for each conclusion of Lemma~\ref{Large-Set}, there
exists a small matching $M'$ in $H$ such that for each legal $k$-set $S$,
there are at least $k$-pairwise disjoint $S$-absorbing matchings in $H$.

\begin{lemma}\label{Absorb-lem}
Let $k \geq 3$ be a positive integer.
There exists $n_2 > 0$ such that the following holds.
Let $H$ be a
$k$-partite $k$-graph with partition
classes $V_1, \ldots, V_k$ such that $|V_i| = n > n_2$ for $i \in [k]$  and $\delta_{k-1}(H)\geq (1/2-1/\log n)n$.
\begin{itemize}
\item[$(i)$] If for all $i\in [k]$ and  $N_i\subseteq V_i$ with
  $|N_i|\geq (1/2-1/\log n)n$, we have
  $e_H(N_1,\ldots,N_k)\geq n^k/\log^3n$,
then there exists a matching $M'$ in $H$ such that $|M'|=O(\log^5 n)$ and
for every legal $k$-set $S \subseteq V(H)$, there are at least $k$
pairwise disjoint $S$-absorbing $k$-matchings in $M'$.

\item [$(ii)$] If there exists $j\in [k]$ such that $|\Lambda_j|\geq n^{k-1}/\log n$,
then there exists a matching $M'$ in $H$  such that $|M'|=O(\log^5 n)$ and
  for every legal $k$-set $S \subseteq V(H)$, there are at least $k$
  pairwise disjoint $S$-absorbing $(k+1)$-matchings in $M'$.
\end{itemize}
\end{lemma}

\pf First, we prove  $(i)$. Suppose for all $i\in [k]$ and  $N_i\subseteq V_i$ with
  $|N_i|\geq (1/2-1/\log n)n$, we have
  $e_H(N_1,\ldots,N_k)\geq n^k/\log^3n$. So we can apply $(i)$ of Lemma \ref{absorb-num}.

For each legal $k$-set $S \subseteq V(H)$,
let $\Gamma(S)$ be the set of  $(S_1,\ldots, S_k)$ with $S_i\subseteq V_i$ and $|S_i|=k$ for $i\in
[k]$ such that $H[\cup_{i\in [k]}S_i]$ has a perfect matching, say $M_{(S_1,\ldots, S_k)}$.
 Then by $(i)$ of Lemma \ref{absorb-num},  $|\Gamma(S)|  =\Omega(n^{k^2}/\log^3 n)/k^k$.
 So there exists $\alpha:=\alpha(k) > 0$ such that $|\Gamma(S)| \geq \alpha {n\choose k}^k/ \log^3 n$.

Let $\mathcal{F}$ be the (random)  family  whose members are  $(S_1,\ldots, S_k)$ with $S_i\subseteq V_i$ and $|S_i|=k$ for $i\in
[k]$, obtained by choosing each of the $\binom{n}{k}^k$ such  $(S_1,\ldots, S_k)$ independently with probability
\[
p=\frac{\log^5 n}{{n\choose k}^k}.
\]
Note that $p<1$ as we can choose $n_2$ large enough. Then
\[
\mathbb{E}(|\mathcal{F}|)=p  {n\choose k}^k=\log^5 n,
\]
and for each legal $k$-set $S\subseteq V(H)$,
\[
\mathbb{E}(|\mathcal{F} \cap \Gamma(S)|) \geq p  \alpha {n\choose k}^k/\log^3 n = \alpha \log^2 n.
\]

By Lemma~\ref{chernoff} and by choosing $n_2$ large enough, we have, for $n>n_2$,
$$
\mathbb{P}[|\mathcal{F}| > 2 \log^5 n] = \mathbb{P}[|\mathcal{F}| > 2 \mathbb{E}(|\mathcal{F}|)] \leq e^{-\mathbb{E}(|\mathcal{F}|)/3} = e^{- (\log^5 n)/3} < 1/10.
$$
So with probability at least $9/10$
\begin{align}\label{prop-F}
|\mathcal{F}|\leq 2 \log^5 n.
\end{align}
Again by Lemma~\ref{chernoff} and  by choosing $n_2$ large enough, we have, for $n>n_2$,
\begin{align*}
\mathbb{P}[|\mathcal{F} \cap \Gamma(S)| \leq (\alpha \log^2 n)/2 ] &\leq  \mathbb{P}[|\mathcal{F} \cap \Gamma(S)| \leq  \mathbb{E}(|\mathcal{F} \cap \Gamma(S)|)/2 ] \\
& \leq  e^{-\mathbb{E}(|\mathcal{F} \cap \Gamma(S)|)/8}\\
&\leq  e^{- (\alpha\log^2 n)/8}.
\end{align*}
So by union bound and choosing $n_2$  large, we have for $n>n_2$,
$$
\mathbb{P}[\exists \text{ legal } S \subseteq V(H): |\mathcal{F} \cap \Gamma(S)| \leq (\alpha \log^2 n)/2 ] \leq n^k e^{- (\alpha\log^2 n)/8} = 2n^{k - (\alpha\log n)/8} < 1/10.
$$
Thus,  with probability at least $9/10$, for each legal $k$-set
$S\subseteq V(H)$, we have
\begin{align}\label{prop-F-S}
|\mathcal{F} \cap \Gamma(S)| \geq (\alpha\log^2 n)/2 > k.
\end{align}

Furthermore, the expected number of pairs of elements $(S_1, \ldots, S_k), (T_1, \ldots, T_k)\in {\cal F}$
satisfying $(\cup_{i \in [k]} S_i) \cap (\cup_{i \in [k]} T_i) \neq \emptyset$ is at most
\[
{n\choose k}^k k^2 {n-1 \choose k-1}  {n \choose k}^{k-1}  p^2 \leq \frac{k^3\log ^{10} n}{n}<1/2.
\]
Thus, with probability at least $1/2$ (by Markov's inequality), for all distinct $(S_1,\ldots, S_k)\in {\cal F}$ and
$(T_1,\ldots, T_k)\in \mathcal{F}$,
\begin{align}\label{F-cap}
\mbox{$\cup_{i\in[k]} S_i$ and $\cup_{i\in [k]} T_i$ are disjoint.}
\end{align}

Hence, with positive probability, $\mathcal{F}$ satisfies (\ref{prop-F}),  (\ref{prop-F-S}), and (\ref{F-cap}). So we may assume that ${\cal F}$
satisfies (\ref{prop-F}),  (\ref{prop-F-S}), and (\ref{F-cap}).
Let $M'$ be the union of $M_{(S_1,\ldots, S_k)}$ for $(S_1, \ldots, S_k)\in {\cal F}$.
Then $M'$ is the desired matching for $(i)$.

\medskip

Next we prove $(ii)$.  Suppose there exists $j\in [k]$ such that
$|\Lambda_j|\geq n^{k-1}/\log n$; so that we can apply $(ii)$ of Lemma
\ref{absorb-num}.

For each legal $k$-set $S \subseteq V(H)$,
let $\Gamma'(S)$ be the set of  sequences $(S_1,\ldots, S_k)$, with
$S_i\subseteq V_i$ and $|S_i|=k+1$ for $i\in [k]$, such that $H[\cup_{i\in [k]}S_i]$ has a perfect matching, say $M_{(S_1,\ldots, S_k)}$.
 Then by  $(ii)$ of Lemma \ref{absorb-num},  $|\Gamma'(S)|  =\Omega(n^{k^2+k}/\log^3 n)/(k+1)^k$.
 So there exists $\alpha' > 0$ such that $|\Gamma'(S)| \geq \alpha'  {n\choose
 {k+1}}^k/ \log^3 n$.

We form a random family $\mathcal{G}$ consisting of sequences $(S_1,\ldots, S_k)$, with
$S_i\subseteq V_i$ and $|S_i|=k+1$ for $i\in [k]$, by selecting each
of the ${n\choose k+1}^k$ such
$(S_1,\ldots, S_k)$ independently with probability
\[
p=\frac{\log^5 n}{{n\choose k+1}^k}.
\]
Note that $p<1$ by choosing $n_2$ large enough. Then
\[
\mathbb{E}(|\mathcal{G}|)=p {n\choose k+1}^k=\log^5 n,
\]
and for each legal $k$-set $S\subseteq V(H)$,
\[
\mathbb{E}(|\mathcal{G} \cap \Gamma'(S)|) \geq p \alpha'  {n\choose
 {k+1}}^k/\log^3 n=\alpha' \log^2 n.
\]
By Lemma~\ref{chernoff} and by choosing $n_2$ large enough, we have for $n>n_2$,
$$
\mathbb{P}[|\mathcal{G}| > 2 \log^5 n] = \mathbb{P}[|\mathcal{G}| > 2 \mathbb{E}(|\mathcal{G}|)] \leq 2e^{-\mathbb{E}(|\mathcal{G}|)/3} =2 e^{- (\log^5 n)/3} < 1/10.
$$
So with probability at least $9/10$,
\begin{align}\label{prop-G}
|\mathcal{G}|\leq 2 \log^5 n.
\end{align}
Again by Lemma~\ref{chernoff} and by choosing $n_2$ large enough, we have for $n>n_2$,
\begin{align*}
\mathbb{P}[|\mathcal{G} \cap \Gamma'(S)| \leq (\alpha' \log^2 n)/2 ] & \leq  \mathbb{P}[|\mathcal{G} \cap \Gamma'(S)|
\leq \mathbb{E}(|\mathcal{G} \cap \Gamma'(S)|/2)]\\
&\leq  e^{-\mathbb{E}(|\mathcal{G} \cap \Gamma'(S)|)/8} \\
& \leq  e^{- \alpha'\log^2 n/8}.
\end{align*}
So  by union bound and choosing $n_2$ large,
$$
\mathbb{P}[\exists \text{ legal } S \subseteq V(H): | \mathcal{G} \cap \Gamma'(S)| \leq \alpha' \log^2 n/2] \leq n^k e^{- (\alpha'\log^2 n)/8} = n^{k - (\alpha'\log n)/8} < 1/10.
$$
Hence,  with probability at least $9/10$, for each legal $k$-set
$S\subseteq V(H)$,
\begin{align}\label{prop-G-S}
|\mathcal{G} \cap \Gamma'(S)| \geq (\alpha' \log^2 n)/2 > k.
\end{align}

Furthermore,  the expected number of pairs $(S_1,\ldots,
S_k),(T_1,\ldots, T_k)\in \mathcal{G}$ with $(\cup_{i\in [k]} S_i)\cap
(\cup_{i\in [k]}T_i)\ne \emptyset$ is
\[
{n\choose k+1}^k  k(k+1)  {n-1 \choose k}  {n\choose k+1}^{k-1} p^2 \leq \frac{(k+1)^3\log ^{10} n}{n}< 1/2.
\]
Thus, by Markov's inequality, with probability at least $1/2$, for all distinct $(S_1,\ldots,
S_k)\in {\cal G}$ and $(T_1,\ldots, T_k)\in \mathcal{G}$,
\begin{align}\label{G-cap}
\mbox{ $(\cup_{i\in [k]} S_i)\cap
(\cup_{i\in [k]}T_i)= \emptyset$.}
\end{align}

Hence, with positive probability, $\mathcal{G}$ satisfies (\ref{prop-G}),  (\ref{prop-G-S}), and (\ref{G-cap}).
So we may assume that $\mathcal{G}$ satisfies (\ref{prop-G}),  (\ref{prop-G-S}), and (\ref{G-cap})
Let $M'$ be the union of $M'_{(S_1,\ldots, S_k)}$ for all $(S_1,\ldots, S_k)\in {\cal G}$.
Now $M'$ gives the desired matching for $(ii)$. \qed

\medskip

\begin{coro}
\label{lemma_away_from_extremal}
Let $k \geq 3$ be a positive integer.
For any $\varepsilon>0$, there exists $n_3 > 0$ such that the following holds.
Let $H$ be a $k$-partite $k$-graph with $n>n_3$ vertices in each partition class.
Suppose $\delta_{k-1}(H)\geq (1/2-1/\log n )n$ and $H$ is not $\varepsilon$-close to $H_0(k,n)$.
Then $H$ has a perfect matching.
\end{coro}

\pf Choose $n_3$ large enough so that we can apply Lemmas  \ref{near_perfect}, \ref{Large-Set}, and \ref{Absorb-lem}.

By Lemmas \ref{Large-Set} and \ref{Absorb-lem}, $H$ contains a matching $M$ such that $|M| \leq \beta\log^5 n$ for some constant $\beta > 0$
(dependent on $k$ only) and, for every
legal $k$-set $S \subseteq V(H)$, there are at least $k$ disjoint $S$-absorbing $k$-matchings in $M$, or for every
legal $k$-set $S \subseteq V(H)$, there are at least $k$ disjoint $S$-absorbing $(k+1)$-matchings in $M$.

For $k \geq 3$,
\begin{align*}
\delta_{k-1}(H-V(M))\geq (1/2 - 1/\log n)n - \beta\log^5 n
> n/k,
\end{align*}
where the last inequality holds for $n>n_3$ by choosing $n_3$ large enough.
Thus by Lemma \ref{near_perfect}, $H-V(M)$ contains a  matching $M'$ of size at least $n-|M| - 1$.
Let $S := H - V(M \cup M')$.
If $S = \emptyset$, then $M \cup M'$ is a perfect matching in $H$.
So assume that $S\ne \emptyset$; then $S$ is a legal $k$-set.
Hence
$H[S\cup V(M)]$ has a perfect matching $M''$.
Now $M'\cup M''$ is a perfect matching in $H$.  \qed

\section{Conclusion}

\noindent {\bf Proof of Theorem \ref{main}.}
First, suppose $(i)$ or $(ii)$ holds. Then there exist integers $d_1,\ldots, d_k$ such that $\sum_{i=1}^k d_i$ is odd and $H\subseteq H_0(d_1,d_2,\ldots,d_k;k,n)$.
By definition of $H_0(d_1,d_2,\ldots,d_k;k,n)$,
there exists $D \subseteq V(H)$ such that
$|D| = \sum_{i=1}^k d_i$ is odd and
$|e\cap D|$ is even for all $e\in E(H)$.
Hence, $H$ contains no perfect matchings.

Now assume that $H$ has no perfect matching.
Fix $\varepsilon>0$ so that $$\sqrt{\varepsilon}< \min\{1/(100k^2),1/(k(10k^2)^{k-1})\}.$$
Then by Corollary \ref{lemma_away_from_extremal},  $H$ must be $\varepsilon$-close to $H_0(k,n)$.
Hence by Corollary \ref{near_critical3}, $(i)$ or $(ii)$ holds.
\qed

\end{document}